\documentclass[11pt,twoside]{article} 

\RequirePackage{amsthm,amsmath,amsfonts}

\RequirePackage[colorlinks]{hyperref}

\RequirePackage{hypernat}

\RequirePackage{epsfig, graphicx, color}

\RequirePackage{latexsym}

\RequirePackage{psfrag}

\usepackage{fullpage}

\hypersetup{linkcolor=magenta}

\hypersetup{citecolor=red}

\theoremstyle{plain}

\newtheorem{theo}{Theorem}[section]

\newtheorem{lem}{Lemma}[section]
\newtheorem{prop}{Proposition}[section]
\newtheorem{cor}{Corollary}[section]

\theoremstyle{definition} 

\newtheorem{nota}{Notation}[section]
\newtheorem{de}{Definition}[section]
\newtheorem{exa}{Example}[section]
\newtheorem{as}{Assumption}[section]
\newtheorem{alg}{Algorithm}[section]

\newcommand{\btheo}{\begin{theo}}
\newcommand{\bde}{\begin{de}}
\newcommand{\ble}{\begin{lem}}
\newcommand{\bpr}{\begin{prop}}
\newcommand{\bno}{\begin{nota}}
\newcommand{\bex}{\begin{exa}}
\newcommand{\bcor}{\begin{cor}}
\newcommand{\spro}{\begin{proof}}
\newcommand{\bas}{\begin{as}}
\newcommand{\balg}{\begin{alg}}

\newcommand{\etheo}{\end{theo}}
\newcommand{\ede}{\end{de}}
\newcommand{\ele}{\end{lem}}
\newcommand{\epr}{\end{prop}}
\newcommand{\eno}{\end{nota}}
\newcommand{\eex}{\end{exa}}
\newcommand{\ecor}{\end{cor}}
\newcommand{\fpro}{\end{proof}}
\newcommand{\eas}{\end{as}}
\newcommand{\ealg}{\end{alg}}

\theoremstyle{plain}

\newtheorem{theos}{Theorem}
\newtheorem{props}{Proposition}
\newtheorem{lems}{Lemma}
\newtheorem{cors}{Corollary}

\theoremstyle{definition}
\newtheorem{exas}{Example}
\newtheorem{algs}{Algorithm}
\newtheorem{asss}{Assumption}
\newtheorem{defns}{Definition}

\newcommand{\btheos}{\begin{theos}}
\newcommand{\etheos}{\end{theos}}
\newcommand{\bprops}{\begin{props}}
\newcommand{\eprops}{\end{props}}
\newcommand{\bdes}{\begin{defns}}
\newcommand{\edes}{\end{defns}}
\newcommand{\blems}{\begin{lems}}
\newcommand{\elems}{\end{lems}}
\newcommand{\bcors}{\begin{cors}}
\newcommand{\ecors}{\end{cors}}
\newcommand{\bexs}{\begin{exas}}
\newcommand{\eexs}{\end{exas}}
\newcommand{\balgs}{\begin{algs}}
\newcommand{\ealgs}{\end{algs}}
\newcommand{\bass}{\begin{asss}}
\newcommand{\eass}{\end{asss}}

\newcommand{\y}{\ensuremath{y}}
\newcommand{\w}{\ensuremath{w}}
\newcommand{\yi}{\ensuremath{\y_i}}
\newcommand{\exi}{\ensuremath{x_i}}

\newcommand{\wi}{\ensuremath{w_i}}
\newcommand{\s}{{\ensuremath{s}}}

\newcommand{\PackEnt}{\ensuremath{\log M}}

\newcommand{\numobs}{\ensuremath{n}}
\newcommand{\pdim}{\ensuremath{d}}
\newcommand{\defn}{\ensuremath{:\,=}}

\newcommand{\epspack}{\ensuremath{\delta_{\numobs}}}

\newcommand{\mP}{\ensuremath{\mathbb{P}}}

\newcommand{\mE}{\ensuremath{\mathbb{E}}}

\newcommand{\Hilb}{\ensuremath{{\mathcal{H}}}}

\newcommand{\f}{\ensuremath{f}}

\newcommand{\fTr}{\ensuremath{f^{*}}}
\newcommand{\fHat}{\ensuremath{\widehat{f}}}
\newcommand{\DeltaHat}{\ensuremath{\widehat{\Delta}}}

\newcommand{\xij}{\ensuremath{x_{ij}}}

\newcommand{\LTP}{\ensuremath{L^2(\mathbb{P})}}
\newcommand{\LTQ}{\ensuremath{L^2(\mathbb{Q})}}
\newcommand{\LTPn}{\ensuremath{{L^2(\mathbb{P}_n)}}}
\newcommand{\LTQn}{\ensuremath{{L^2(\mathbb{Q}_n)}}}

\newcommand{\PackNum}{\ensuremath{M}}

\newcommand{\kull}[2]{\ensuremath{D(#1 \, \| \, #2)}}
\newcommand{\GenSet}{\ensuremath{\mathcal{G}}}

\newcommand{\real}{\ensuremath{\mathbb{R}}}
\newcommand{\Sset}{\ensuremath{S}}

\newcommand{\Hil}{\ensuremath{\mathcal{H}}}
\newcommand{\Xspace}{\ensuremath{\mathcal{X}}}
\newcommand{\mprob}{\ensuremath{\mathbb{P}}}
\newcommand{\Qprob}{\ensuremath{\mathbb{Q}}}

\newcommand{\fstar}{\ensuremath{f^*}}

\newcommand{\order}{\ensuremath{\mathcal{O}}}

\long\def\comment#1{}

\newcommand{\Fclass}{\ensuremath{\mathcal{F}}}

\newcommand{\plaincon}{\ensuremath{c}}

\newcommand{\RegPar}{\ensuremath{\lambda_n}}

\newcommand{\Event}{\mathcal{E}}

\newcommand{\Ind}{\ensuremath{\mathbb{I}}}

\newcommand{\Ball}{\ensuremath{\mathbb{B}}}

\newcommand{\NonPar}{\ensuremath{\alpha}}
\newcommand{\Kern}{\ensuremath{K}}
\newcommand{\ParDim}{\ensuremath{m}}

\newcommand{\Exs}{\ensuremath{\mathbb{E}}}
\newcommand{\epscrit}{\ensuremath{\epsilon_\numobs}}

\newcommand{\AuxEvent}{\ensuremath{\mathcal{A}}}

\newcommand{\xsam}[1]{\ensuremath{x_{#1}}}

\newcommand{\gtil}{\ensuremath{\widetilde{g}}}

\newcommand{\fhat}{\ensuremath{\widehat{f}}}

\newcommand{\var}{\ensuremath{\operatorname{var}}}

\newcommand{\kdim}{\ensuremath{\spindex}}
\newcommand{\spindex}{\ensuremath{s}}

\newcommand{\UniCon}{\ensuremath{C}}

\newcommand{\Ker}{\ensuremath{\mathbb{K}}}

\newcommand{\kerrank}{\ensuremath{m}}

\newcommand{\Xcal}{\ensuremath{\mathcal{X}}}

\newcommand{\LocGauss}{\ensuremath{\mathcal{Q}_{w,\numobs}}}
\newcommand{\LocRad}{\ensuremath{\mathcal{Q}_{\sigma,\numobs}}}
\newcommand{\LocGaussEmp}{\ensuremath{\widehat{\mathcal{Q}}_{w,\numobs}}}

\newcommand{\capxsam}[1]{\ensuremath{x_{#1}}}

\newcommand{\radepl}{\ensuremath{\sigma}}
\newcommand{\rade}[1]{\ensuremath{\radepl_{#1}}}

\newcommand{\regpar}{\ensuremath{\lambda_\numobs}}

\newcommand{\delunihat}[1]{\ensuremath{\widehat{\unicritpl}_{\numobs,
      #1}}}
\newcommand{\delhatuni}[1]{\ensuremath{\delunihat{#1}}}

\newcommand{\regpartwo}{\ensuremath{\rho_\numobs}}
\newcommand{\Sbar}{\ensuremath{{S^c}}}

\newcommand{\gamcrit}{\ensuremath{\gamma_\numobs}}

\newcommand{\Bevent}{\ensuremath{\mathcal{B}}}

\newcommand{\Qrad}{\ensuremath{\mathcal{Q}_\numobs}}

\newcommand{\Rvar}{\ensuremath{Z_\numobs}}
\newcommand{\RvarHat}{\ensuremath{\widehat{Z}_\numobs}}

\newcommand{\unicrit}{\ensuremath{\nu_\numobs}}
\newcommand{\unicritpl}{\ensuremath{\nu}}

\newcommand{\IndSet}{\ensuremath{\mathcal{I}}}
\newcommand{\Hset}{\ensuremath{\mathfrak{S}}}
\newcommand{\Aset}{\ensuremath{\mathcal{A}}}

\newcommand{\mdim}{\ensuremath{m}}
\newcommand{\smooth}{\ensuremath{\alpha}}

\newcommand{\MiniMax}{\ensuremath{\mathfrak{M}}}
\newcommand{\PackVar}{\ensuremath{\Theta}}
\newcommand{\PackVarHat}{\ensuremath{\widehat{\PackVar}}}
\newcommand{\Xdes}{\ensuremath{{X_1^\numobs}}}
\newcommand{\Ysam}{\ensuremath{Y_1^\numobs}}

\newcommand{\MyHclass}{\ensuremath{\Hil(\Sset)}}

\newcommand{\CostFun}{\mathcal{L}}
\newcommand{\TailEvent}{\ensuremath{\mathcal{T}}}

\newcommand{\lpe}{\ensuremath{L^2(\Qprob_\numobs)}}

\newcommand{\Lpeone}[1]{\ensuremath{\|#1\|_{\numobs,1}}}

\newcommand{\Hilone}[1]{\ensuremath{\|#1\|_{\Hil,1}}}

\newcommand{\MiniMaxPop}{\ensuremath{\MiniMax_\mprob(\MyBigClass)}}
\newcommand{\MiniMaxFast}{\ensuremath{\MiniMax_\mprob(\MyNewClass)}}

\newcommand{\inprod}[2]{\ensuremath{\langle #1 , \, #2 \rangle}}

\newcommand{\InfoX}{\ensuremath{I_\Xdes}}

\newcommand{\MyBigClass}{\ensuremath{\Fclass_{\pdim, \spindex, \Hil}}}

\newcommand{\alhat}{\ensuremath{\widehat{\alpha}}}

\newcommand{\SpecM}{\ensuremath{N}}

\newcommand{\GoodPackSize}{M}
\newcommand{\htil}{\ensuremath{\widetilde{h}}}

\newcommand{\AnnoyEvent}{\ensuremath{\mathcal{D}}}

\newcommand{\LossTil}{\ensuremath{\widetilde{\mathcal{L}}}}

\newcommand{\MyNewClass}{\mathcal{F}^*_{\pdim,\s, \Hil}(B)}
\newcommand{\MyNewSmallClass}{\Hil(S,B)}
\newcommand{\TRUNCLASS}{\Hil(S,B, M)}

\newcommand{\Gclass}{\ensuremath{\mathcal{G}}}

\newcommand{\GLOB}{\ensuremath{B}}

\newcommand{\UnifBoundRate}{\delta_n}

\newcommand{\DelHat}{\widehat{\Delta}}

\newcommand{\Nat}{\ensuremath{\mathbb{N}}}

\newcommand{\tracer}[2]{\ensuremath{\langle \!\langle {#1}, \; {#2}
\rangle \!\rangle}}

\newcommand{\HACK}{\ensuremath{\mathcal{H}(S, 2 \GLOB)}}

\newcommand{\Aone}{\ensuremath{A_1}}
\newcommand{\Atwo}{\ensuremath{A_2}} 

\newcommand{\smallbou}{\ensuremath{b}}

\newcommand{\FUNCLASS}{\ensuremath{\mathcal{G}}} 

\newcommand{\tildelcrit}{\tilde{\delta}_n}
\newcommand{\ybar}{\bar{y}_n}

\newcommand{\Nprop}{\ensuremath{N_{\tiny{\mbox{pr}}}}}

\newcommand{\gvar}{\ensuremath{w}}
\newcommand{\newgvar}{\ensuremath{\varepsilon}}

\begin{document}

\begin{center}

{\bf{\LARGE{ Minimax-optimal rates for sparse additive models over
      kernel classes via convex programming}}} \\

\vspace*{.3in}

\begin{tabular}{cc}
Garvesh Raskutti$^1$ &  Martin J. Wainwright$^{1,2}$ \\
	    \texttt{garveshr@stat.berkeley.edu} &
	    \texttt{wainwrig@stat.berkeley.edu}
\end{tabular}

\vspace*{.1in}

\begin{tabular}{c}
Bin Yu$^{1,2}$ \\
\texttt{binyu@stat.berkeley.edu}
\end{tabular}

\vspace*{.2in}

\begin{tabular}{c}
Departments of Statistics$^1$, and EECS$^2$ \\
UC Berkeley,  Berkeley, CA 
\end{tabular}

\vspace*{.2in}

\begin{center}
Revised version of August 2010 paper
\end{center}

\vspace*{.2in}

\begin{abstract}
Sparse additive models are families of $d$-variate functions that have
the additive decomposition \mbox{$f^* = \sum_{j \in S} f^*_j$,} where
$S$ is an unknown subset of cardinality $s \ll d$. In this paper, we
consider the case where each univariate component function $f^*_j$
lies in a reproducing kernel Hilbert space (RKHS), and analyze a
method for estimating the unknown function $f^*$ based on kernels
combined with $\ell_1$-type convex regularization.  Working within a
high-dimensional framework that allows both the dimension $d$ and
sparsity $s$ to increase with $n$, we derive convergence rates (upper
bounds) in the $L^2(\mathbb{P})$ and $L^2(\mathbb{P}_n)$ norms over
the class $\MyBigClass$ of sparse additive models with each univariate
function $f^*_j$ in the unit ball of a univariate RKHS with bounded
kernel function.  We complement our upper bounds by deriving minimax
lower bounds on the $L^2(\mathbb{P})$ error, thereby showing the
optimality of our method.  Thus, we obtain optimal minimax rates for
many interesting classes of sparse additive models, including
polynomials, splines, and Sobolev classes.  We also show that if, in
contrast to our univariate conditions, the multivariate function class
is assumed to be globally bounded, then much faster estimation rates
are possible for any sparsity $s = \Omega(\sqrt{n})$, showing that
global boundedness is a significant restriction in the
high-dimensional setting.
\end{abstract}

\end{center}


\section{Introduction}

The past decade has witnessed a flurry of research on sparsity
constraints in statistical models.  Sparsity is an attractive
assumption for both practical and theoretical reasons: it leads to
more interpretable models, reduces computational cost, and allows for
model identifiability even under high-dimensional scaling, where the
dimension $\pdim$ exceeds the sample size $\numobs$.  While a large
body of work has focused on sparse linear models, many applications
call for the additional flexibility provided by non-parametric models.
In the general setting, a non-parametric regression model takes the
form $\y = \fTr(x_1, \ldots, x_\pdim) + w$, where $f^*:\real^\pdim
\rightarrow \real$ is the unknown regression function, and $w$ is
scalar observation noise.  Unfortunately, this general non-parametric
model is known to suffer severely from the so-called ``curse of
dimensionality'', in that for most natural function classes (e.g.,
twice differentiable functions), the sample size $\numobs$ required to
achieve any given error grows exponentially in the dimension $\pdim$.
Given this curse of dimensionality, it is essential to further constrain
the complexity of possible functions $\fTr$.  One attractive candidate
is the class of \emph{additive non-parametric
  models}~\cite{HasTib86}, in which the function $\fTr$ has an
additive decomposition of the form
\begin{align}
\label{EqnAddModel}
\fTr(x_1,x_2, \ldots, x_{\pdim}) & = \sum_{j=1}^{\pdim}{f_j^{*}(x_j)},
\end{align}
where each component function $f_j^{*}$ is univariate.  Given this
additive form, this function class no longer suffers from the
exponential explosion in sample size of the general non-parametric
model.  Nonetheless, one still requires a sample size $\numobs \gg
\pdim$ for consistent estimation; note that this is true even for the
linear model, which is a special case of equation~\eqref{EqnAddModel}.

A natural extension of sparse linear models is the class of \emph{sparse additive models}, in which the unknown regression function is assumed to have a
decomposition of the form
\begin{align}
\label{EqnSPAM}
\fTr(x_1,x_2 \ldots, x_{\pdim}) & = \sum_{j \in \Sset}{f_j^{*}(x_j)},
\end{align}
where $\Sset \subseteq \{ 1, 2, \ldots, \pdim \}$ is some unknown
subset of cardinality $|S| = \s$.  Of primary interest is the case
when the decomposition is genuinely sparse, so that $\s \ll \pdim$.
To the best of our knowledge, this model class was first introduced in
Lin and Zhang~\cite{LinZhang06}, and has since been studied by various
researchers (e.g.,~\cite{KolYua10Journal,Meier09,RavLiuLafWas08,Yuan07}).
Note that the sparse additive model~\eqref{EqnSPAM} is a natural
generalization of the sparse linear model, to which it reduces when
each univariate function is constrained to be linear.

In past work, several groups have proposed computationally efficient
methods for estimating sparse additive models~\eqref{EqnSPAM}.  Just
as $\ell_1$-based relaxations such as the Lasso have desirable
properties for sparse parametric models, more general $\ell_1$-based
approaches have proven to be successful in this setting. Lin and
Zhang~\cite{LinZhang06} proposed the COSSO method, which extends the
Lasso to cases where the component functions $f_j^*$ lie in a
reproducing kernel Hilbert space (RKHS); see also Yuan~\cite{Yuan07}
for a similar extension of the non-negative garrote~\cite{Breiman95}.
Bach~\cite{Bach08a} analyzes a closely related method for the RKHS
setting, in which least-squares loss is penalized by an $\ell_1$-sum
of Hilbert norms, and establishes consistency results in the classical
(fixed $\pdim$) setting.  Other related $\ell_1$-based methods have
been proposed in independent work by Koltchinskii and
Yuan~\cite{KolYua08}, Ravikumar et al.~\cite{RavLiuLafWas08} and Meier
et al.~\cite{Meier09}, and analyzed under high-dimensional scaling ($n
\ll \pdim$).  As we describe in more detail in
Section~\ref{SecCompare}, each of the above papers establish
consistency and convergence rates for the prediction error under
certain conditions on the covariates as well as the sparsity $\s$ and
dimension $\pdim$.  However, it is not clear whether the rates
obtained in these papers are sharp for the given methods, nor whether
the rates are minimax-optimal.  Past work by Koltchinskii and
Yuan~\cite{KolYua10Journal} establishes rates for sparse additive
models with an additional global boundedness condition, but as will be
discussed at more length in the sequel, these rates are not minimax
optimal in general.

This paper makes three main contributions to this line of research.
Our first contribution is to analyze a simple polynomial-time method
for estimating sparse additive models and provide upper bounds on the
error in the $\LTP$ and $\LTPn$ norms. The estimator\footnote{The
  estimator is the same as the estimator analyzed in Koltchinskii and
  Yuan~\cite{KolYua10Journal}. We proposed the estimator concurrently
  (see our earlier preprint~\cite{RasWaiYu10}) and as we discuss later
  analyze the same estimator under less restrictive conditions than
  those imposed in Koltchinskii and Yuan~\cite{KolYua10Journal}.} we
analyze is based on a combination of least-squares loss with two
$\ell_1$-based sparsity penalty terms, one corresponding to an
$\ell_1/\LTPn$ norm and the other an $\ell_1/\|\cdot\|_\Hil$ norm. Our
first main result (Theorem~\ref{ThmPolyAchievable}) shows that with
high probability, if we assume the univariate functions are bounded
and independent, the error of our procedure in the squared $\LTPn$ and
$\LTP$ norms is bounded by $\order\big(\frac{\s \log \pdim}{\numobs}+
\s \unicrit^2\big)$, where the quantity $\unicrit^2$ corresponds to
the optimal rate for estimating a single univariate
function. Importantly, our analysis does \emph{not} require a global
boundedness condition on the class $\MyBigClass$ of all
$\spindex$-sparse models, an assumption that is often imposed in
classical non-parametric analysis.  Indeed, as we discuss below, when
such a condition is imposed, then significantly faster rates of
estimation are possible.  The proof of
Theorem~\ref{ThmPolyAchievable} involves a combination of techniques
for analyzing $M$-estimators with decomposable
regularizers~\cite{Neg09} and techniques in empirical process theory
for analyzing kernel classes~\cite{Bartlett05,Mendelson02,vandeGeer}.

Our second contribution is complementary in nature, in that it
establishes algorithm-independent minimax lower bounds on $\LTP$
error. These minimax lower bounds, stated in Theorem~\ref{ThmLower},
are specified in terms of the metric entropy of the underlying
univariate function classes. For both finite-rank kernel classes and
Sobolev-type classes, these lower bounds match our achievable results,
as stated in Corollaries~\ref{CorAchieveFinite}
and~\ref{CorAchieveSmooth}, up to constant factors in the regime of
sub-linear sparsity ($\spindex = o(\pdim)$).  Thus, for these function
classes, we have a sharp characterization of the associated minimax
rates. The lower bounds derived in this paper initially appeared in
the Proceedings of the NIPS Conference (December 2009). The proofs of
Theorem 2 is based on characterizing the packing entropies of the
class of sparse additive models, combined with classical information
theoretic techniques involving Fano's inequality and variants (see,
e.g. the papers~\cite{Hasminskii78,YanBar99,Yu}).

Our third contribution is to determine upper bounds on minimax $\LTP$
and $\LTPn$ error when we impose a global boundedness assumption on
the class $\MyBigClass$, meaning that the quantity
\mbox{$\GLOB(\MyBigClass) = \sup_{f \in \MyBigClass} \sup_x
  |\sum_{j=1}^\pdim f_j(x_j)|$} is assumed to be bounded independently
of $(\s, \pdim)$.  As mentioned earlier, our upper bound in
Theorem~\ref{ThmPolyAchievable} does \emph{not} impose a global
boundedness condition, whereas in contrast, past work by Koltchinskii
and Yuan~\cite{KolYua10Journal} did impose such a global boundedness
condition in their analysis of the same $\ell_1$-kernel-based
estimator.  Under global boundedness, their work provides rates on the
$L^2(\mathbb{P})$ and $L^2(\mathbb{P}_n)$ norm that are of the same
order as the results presented here.  It is natural to wonder whether
or not this difference is actually significant---that is, do the
minimax rates for the class of sparse additive models depend on
whether or not global boundedness is imposed?  In
Section~\ref{SecKYCompare}, we define the class $\MyNewClass$ of
sparse additive models with the additional assumption that
$\GLOB(\MyBigClass) \leq B$.  Theorem~\ref{ThmFast} and
Corollary~\ref{CorFast} in this paper provide upper bounds on the
minimax rate in $\LTP$ and $\LTPn$ error over the class $\MyNewClass$.
These rates are faster than those of Theorem 3 in the
paper~\cite{KolYua10Journal} (KY), showing that the KY rates are not
minimax optimal for problems with $\s = \Omega(\sqrt{n})$.  In this
way, we see that the assumption of global boundedness, though
relatively innocuous for classical (low-dimensional) non-parametric
problems, can be quite limiting in high dimensions.

The remainder of the paper is organized as follows.  In
Section~\ref{SecBackground}, we provide background on kernel spaces
and the class of sparse additive models considered in this paper.
Section~\ref{SecMain} is devoted to the statement of our main results
and discussion of their consequences; it includes description of our
method, the upper bounds on the convergence rate that it achieves, and
a matching set of minimax lower bounds. Section~\ref{SecKYCompare}
emphasizes the restrictiveness of the global uniform boundedness
assumption and in particular Theorem~\ref{ThmFast} and
Corollary~\ref{CorFast} show that there are classes of Sobolev spaces
where under the scaling $\s = \Omega(\sqrt{n})$, optimal rates of
convergence are faster than rates proven in
Theorem~\ref{ThmLower}. Section~\ref{SecProofs} is devoted to the
proofs of our three main theorems, with the more technical details
deferred to the Appendices.  We conclude with a discussion in
Section~\ref{SecDiscuss}.

\section{Background and problem set-up}
\label{SecBackground}

We begin with some background on reproducing kernel Hilbert spaces,
before providing a precise definition of the class of sparse additive
models studied in this paper.

\subsection{Reproducing kernel Hilbert spaces}

Given a subset $\Xspace \subset \real$ and a probability measure
$\Qprob$ on $\Xspace$, we consider a Hilbert space \mbox{$\Hil \subset
\LTQ$,} meaning a family of functions $g: \Xcal \rightarrow \real$,
with $\|g\|_{\LTQ} < \infty$, and an associated inner product
$\inprod{\cdot}{\cdot}_\Hil$ under which $\Hil$ is complete.  The
space $\Hil$ is a reproducing kernel Hilbert space (RKHS) if there
exists a symmetric function $\Ker: \Xcal \times \Xcal \rightarrow
\real_+$ such that: (a) for each $x \in \Xspace$, the function
$\Ker(\cdot, x)$ belongs to the Hilbert space $\Hil$, and (b) we have
the reproducing relation $f(x) = \inprod{f}{\Ker(\cdot, x)}_{\Hil}$
for all $f \in \Hil$.  Any such kernel function must be positive
semidefinite; under suitable regularity conditions, Mercer's
theorem~\cite{Mercer09} guarantees that the kernel has an
eigen-expansion of the form
\begin{align}
\label{EqnMercer}
\Ker(x, x') & = \sum_{k=1}^\infty \mu_k \phi_k(x)
\phi_\ell(x'),
\end{align}
where $\mu_1 \geq \mu_2 \geq \mu_3 \geq \ldots \geq 0$ are a
non-negative sequence of eigenvalues, and $\{\phi_k\}_{k=1}^\infty$
are the associated eigenfunctions, taken to be orthonormal in $\LTQ$.
The decay rate of these eigenvalues will play a crucial role in our
analysis, since they ultimately determine the rate $\unicrit$ for the
univariate RKHS's in our function classes.

Since the eigenfunctions $\{\phi_k\}_{k=1}^\infty$ form an
orthonormal basis, any function $f \in \Hil$ has an expansion of the
$f(x) = \sum_{k=1}^{\infty} a_{k} \phi_k(x)$, where $a_k =
\inprod{f}{\phi_k}_{\LTQ} = \int_\Xspace f(x) \phi_k(x) \, d\,
\Qprob(x)$ are (generalized) Fourier coefficients.  Associated with
any two functions in $\Hil$---say \mbox{$f = \sum_{k=1}^\infty
a_k \phi_k$} and \mbox{$g = \sum_{k=1}^\infty b_k
\phi_k$}---are two distinct inner products.  The first is the usual
inner product in the space $\LTQ$---namely,
\mbox{$\inprod{f}{g}_{\LTQ} \defn \int_\Xspace f(x) g(x) \, d\,
\Qprob(x)$.} By Parseval's theorem, it has an equivalent
representation in terms of the expansion coefficients---namely
\begin{align*}
\inprod{f}{g}_{\LTQ} & = \sum_{k=1}^\infty a_k b_k.
\end{align*}
The second inner product, denoted $\inprod{f}{g}_{\Hil}$, is the one
that defines the Hilbert space; it can be written in terms of the
kernel eigenvalues and generalized Fourier coefficients as
\begin{align*}
\inprod{f}{g}_\Hil & = \sum_{k=1}^\infty \frac{a_k
  b_k}{\mu_k}.
\end{align*}
Using this definition, the Hilbert ball of radius $1$ for the Hilbert space $\Hil$ with eigenvalues $\mu_k$ and eigenfunctions $\phi_k(\cdot)$, is:
\begin{equation}
\label{DefnHilbBall}	
\Ball_\Hil(1) = \{f \in \Hil;\;f(\cdot) = \sum_{k=1}^{\infty} a_{k} \phi_k(\cdot)\;|\;   \sum_{k=1}^{\infty} \frac{a_k^2}{\mu_k} \leq 1\}.	
\end{equation}	
For more background on reproducing kernel Hilbert spaces, we refer
the reader to various standard references~\cite{Aronszajn50,Weinert82,Saitoh88,Scholkopf02,Wahba}.

\subsection{Sparse additive models over RKHS}

For each $j = 1, \ldots, \pdim$, let $\Hil_j \subset \LTQ$ be a
reproducing kernel Hilbert space of univariate functions on the domain
$\Xspace \subset \mathbb{R}$.  We assume that
\begin{equation*}
\mE[f_j(x)] = \int_{\Xspace} f_j(x) d \, \Qprob(x) \; = \; 0 \qquad
\mbox{for all $f_j \in \Hilb_j$, and for each $j = 1, 2, \ldots,
  \pdim$. }
\end{equation*}
As will be clarified momentarily, our observation
model~\eqref{EqnLinObs} allows for the possibility of a non-zero mean
$\mu$, so that there is no loss of generality in this assumption.  For
a given subset $\Sset \subset \{ 1, 2, \ldots, \pdim \}$, we define
\begin{align}
\label{EqnDefnMyHclass}
\MyHclass & \defn \big \{ f = \sum_{j \in \Sset} f_j \, \mid \, f_j
\in \Hil_j, \mbox{ and } f_j \in \Ball_{\Hil_j}(1) \; \, \forall \; j
\in \Sset \big \},
\end{align}
corresponding to the class of functions $f: \Xcal^\pdim \rightarrow
\real$ that decompose as sums of univariate functions on co-ordinates
lying within the set $\Sset$.  Note that $\MyHclass$ is also (a subset
of) a reproducing kernel Hilbert space, in particular with the norm
\begin{align*}
\|f\|^2_{\MyHclass} & = \sum_{j \in \Sset} \|f_j\|_{\Hil_j}^2,
\end{align*}
where $\|\cdot\|_{\Hil_j}$ denotes the norm on the univariate Hilbert
space $\Hil_j$.  Finally, for a cardinality $\kdim \in \{1, 2, \ldots,
\lfloor \pdim/2 \rfloor \}$, we define the function class
\begin{align}
\label{EqnDefnFclass}
\MyBigClass & \defn \bigcup_{\substack{\Sset \subset \{1, 2, \ldots,
    \pdim\} \\ |\Sset| = \kdim}} \MyHclass.
\end{align}
To ease notation, we frequently adopt the shorthand $\Fclass =
\MyBigClass$, but the reader should recall that $\Fclass$ depends on
the choice of Hilbert spaces $\{\Hilb_j\}_{j=1}^\pdim$, and moreover,
that we are actually studying a \emph{sequence of function classes}
indexed by $(\pdim, \kdim)$.

Now let $\mprob = \Qprob^\pdim$ denote the product measure on the
space $\Xspace^\pdim \subseteq \real^\pdim$.  Given an arbitrary
$\fstar \in \Fclass$, we consider the observation model
\begin{align}
\label{EqnLinObs}
\yi & = \mu + \fstar(\exi) + \wi, \quad \mbox{for $i = 1,2, \ldots,
  \numobs$,}
\end{align}
where $\{\wi\}_{i=1}^\numobs$ is an i.i.d. sequence of standard normal
variates, and $\{\exi\}_{i=1}^\numobs$ is a sequence of design points
in $\real^\pdim$, sampled in an i.i.d. manner from $\mprob$.

Given an estimate $\fhat$, our goal is to bound the error $\fhat -
\fstar$ under two norms.  The first is the \emph{usual $\LTP$ norm} on
the space $\Fclass$; given the product structure of $\mprob$ and the
additive nature of any $f \in \Fclass$, it has the additive
decomposition $\|f\|_{\LTP}^2 = \sum_{j=1}^\pdim \|f_j\|_{\LTQ}^2$.
In addition, we consider the error in the \emph{empirical
$\LTPn$-norm} defined by the sample $\{\exi\}_{i=1}^\numobs$, defined
as 
\begin{equation*}
\|f\|_{\LTPn}^2 \defn \frac{1}{\numobs} \sum_{i=1}^\numobs
f^2(\exi). 
\end{equation*} 
Unlike the $\LTP$ norm, this norm does not decouple
across the dimensions, but part of our analysis will establish an
approximate form of such decoupling.  For shorthand, we frequently use
the notation $\|f\|_2 = \|\f\|_{\LTP}$ and $\|f\|_n = \|f\|_{\LTPn}$
for a $\pdim$-variate function $f \in \Fclass$.  With a minor abuse of
notation, for a univariate function $f_j \in \Hil_j$, we also use the
shorthands $\|f_j\|_2 = \|f_j\|_{\LTQ}$ and $\|f_j\|_\numobs =
\|f_j\|_{\LTQn}$.


\section{Main results and their consequences}
\label{SecMain}

This section is devoted to the statement of our three main results,
and discussion of some of their consequences.  We begin in
Section~\ref{SecMethod} by describing a regularized $M$-estimator for
sparse additive models, and we state our upper bounds for this
estimator in Section~\ref{SecAchThm}.  This estimator is essentially
equivalent to that analyzed in the paper KY~\cite{KolYua10Journal},
except that we allow for a non-zero mean for the function, and
estimate it as well.  We illustrate our upper bounds for various
concrete instances of kernel classes. In Section~\ref{SecMiniLower},
we state minimax lower bounds on the $\LTP$ error over the class
$\MyBigClass$, which establish the optimality of our procedure. In
Section~\ref{SecCompare}, we provide a detailed comparison between our
results to past work, and in Section~\ref{SecKYCompare} we discuss the
effect of global boundedness conditions on optimal rates.

\subsection{A regularized $M$-estimator for sparse additive models}
\label{SecMethod}

For any function of the form $f = \sum_{j=1}^\pdim {f_j}$, the $(\lpe,
1)$ and $(\Hil,1)$-norms are given by
\begin{equation}
\label{EqnDefnNorms}
\Lpeone{f} \defn \sum_{j = 1}^\pdim \|f_j\|_\numobs, \quad \mbox{and}
\quad \Hilone{f} \defn \sum_{j=1}^\pdim \|f_j\|_\Hil,
\end{equation}
respectively.  Using this notation and defining the sample mean $\ybar = \frac{1}{n}\sum_{i=1}^{n}{y_i}$, we define the cost functional
\begin{align}
\label{EqnDefnCostFun}
\CostFun(f) & = \frac{1}{2 \numobs} \sum_{i=1}^\numobs \big(\yi -\ybar - f(\exi) \big)^2 + \regpar \|f\|_{\numobs,1} + \regpartwo
\|f\|_{\Hil,1}.
\end{align}
The cost functional $\CostFun(f)$ is least-squares loss with a
sparsity penalty $\|f\|_{\numobs,1}$ and a smoothness penalty
$\|f\|_{\Hil,1}$. Here $(\regpar, \regpartwo)$ are a pair of positive
regularization parameters whose choice will be specified by our
theory.  Given this cost functional, we then consider the
$M$-estimator
\begin{align}
\label{EqnMest}
\fhat & \in \arg \min_{f} \CostFun(f) \quad \mbox{subject to $f =
  \sum_{j=1}^\pdim f_j$ and $\|f_j\|_\Hil \leq 1$ for all $j = 1, 2,
  \ldots, \pdim$.}
\end{align}
In this formulation~\eqref{EqnMest}, the problem is
infinite-dimensional in nature, since it involves optimization over
Hilbert spaces.  However, an attractive feature of this $M$-estimator
is that, as a straightforward consequence of the representer
theorem~\cite{Kimeldorf71,Scholkopf02}, it can be reduced to an
equivalent convex program in $\real^{\numobs} \times \real^\pdim$.  In
particular, for each $j = 1, 2, \ldots, \pdim$, let $\Ker^j$ denote
the kernel function for co-ordinate $j$.  Using the notation $x_i =
(x_{i1}, x_{i2}, \ldots, x_{i \pdim})$ for the $i^{th}$ sample, we
define the collection of empirical kernel matrices $\Kern^j \in
\real^{\numobs \times \numobs}$, with entries $\Kern^j_{i \ell} =
\Ker^j(x_{ij}, x_{\ell j})$.  By the representer theorem, any solution
$\fhat$ to the variational problem~\eqref{EqnMest} can be written in
the form
\begin{align*}
\fhat(z_1, \ldots, z_\pdim) & = \sum_{i=1}^\numobs \sum_{j=1}^\pdim
\alhat_{ij} \Ker^j(z_j, x_{ij}),
\end{align*}
for a collection of weights $\big \{\alhat_j \in \real^\numobs, \; j=
1, \ldots, \pdim \big \}$.  The optimal weights are obtained by
solving the convex program
\begin{align}
\label{EqnMestKer}
(\alhat_1, \ldots, \alhat_\pdim) & \in \arg \min_{ \substack{\alpha_j
    \in \real^\numobs \\ \alpha_j^T \Kern^j \alpha_j \leq 1 }} \biggr
    \{ \frac{1}{2 \numobs} \|\y - \ybar - \sum_{j=1}^{\pdim} \Kern^j \alpha_j
    \|_2^2 + \RegPar \sum_{j=1}^{\pdim} \sqrt{\frac{1}{n}\|\Kern^j
    \alpha_j\|_2^2} + \regpartwo \sum_{j=1}^{\pdim}{\sqrt{\alpha_j^T
    \Kern^j \alpha_j}} \biggr \}.
\end{align}
This problem is a second-order cone program (SOCP), and there are
various algorithms for finding a solution to arbitrary accuracy in
time polynomial in $(\numobs, \pdim)$, among them interior point
methods (e.g., see \S 11 in the book~\cite{Boyd02}). \\

Various combinations of sparsity and smoothness penalties have been
used in past work on sparse additive models.  For instance, the method
of Ravikumar et. al~\cite{RavLiuLafWas08} is based on least-squares
loss regularized with single sparsity constraint, and separate
smoothness constraints for each univariate function. They solve the
resulting optimization problem using a back-fitting
procedure. Koltchinskii and Yuan~\cite{KolYua08} develop a method
based on least-squares loss combined with a single penalty term
$\sum_{j=1}^{\pdim}{\|f_j\|_{\Hilb}}$.  Their method also leads to an
SOCP if $\Hilb$ is a reproducing kernel Hilbert space, but differs
from the program~\eqref{EqnMestKer} in lacking the additional sparsity
penalties.  Meier et. al~\cite{Meier09} analyzed least-squares
regularized with a penalty term of the form
$\sum_{j=1}^{\pdim}{\sqrt{\lambda_1\|f_j\|_n^2+ \lambda_2
    \|f_j\|_{\Hilb}^2}}$, where $\lambda_1$ and $\lambda_2$ are a pair
of regularization parameters. In their method, $\lambda_1$ controls
the sparsity while $\lambda_2$ controls the smoothness. If $\Hilb$ is
an RKHS, the method in Meier et. al~\cite{Meier09} reduces to an
ordinary group Lasso problem on a different set of variables, which
can be cast as a quadratic program. The more recent work of
Koltchinskii and Yuan ~\cite{KolYua10Journal} is based on essentially
the same estimator as problem~\eqref{EqnMest}, but they impose
stronger assumptions in their analysis. We provide a more in-depth
comparison of our analysis and results with the past work listed above
in Sections~\ref{SecCompare} and~\ref{SecKYCompare}.

\subsection{Upper bound}
\label{SecAchThm}

We now state a result that provides upper bounds on the estimation
error achieved by the estimator~\eqref{EqnMest}, or
equivalently~\eqref{EqnMestKer}.  To simplify presentation, we state
our result in the special case that the univariate Hilbert space
$\Hil_j, j = 1, \ldots, \pdim$ are all identical, denoted by $\Hil$.
However, the analysis and results extend in a straightforward manner
to the general setting of distinct univariate Hilbert spaces, as we
discuss following the statement of Theorem~\ref{ThmPolyAchievable}.

Let $\mu_1 \geq \mu_2 \geq \ldots \geq 0$ denote the non-negative
eigenvalues of the kernel operator defining the univariate Hilbert
space $\Hilb$, as defined in equation~\eqref{EqnMercer}, and define
the function
\begin{align}
\label{EqnDefnUni}
\LocRad(t) & \defn \frac{1}{\sqrt{\numobs}} \big[ \sum_{\ell=1}^\infty
\min\{t^2, \mu_\ell \} \big]^{1/2}.
\end{align}
Let $\unicrit > 0$ be the smallest positive solution to the inequality
\begin{equation}
\label{EqnCritCon}
40 \unicrit^2 \geq \LocRad(\unicrit),
\end{equation}
where the $40$ is simply used for technical convenience. We refer to
$\unicrit$ as the \emph{critical univariate rate}, as it is the
minimax-optimal rate for $\LTP$-estimation of a single univariate
function in the Hilbert space $\Hilb$ (e.g., ~\cite{Mendelson02,
  vandeGeer}). This quantity will be referred to throughout the
remainder of the paper.

Our choices of regularization parameters are specified in terms of the
quantity
\begin{align}
\label{EqnGamCrit}
\gamcrit & \defn \kappa \max \big \{ \unicrit, \sqrt{\frac{\log
    \pdim}{\numobs}} \big \},
\end{align}
where $\kappa$ is a fixed constant that we choose later. We assume
that each function within the unit ball of the univariate Hilbert
space is uniformly bounded by a constant multiple of its Hilbert
norm---that is, for each $j = 1, \ldots, \pdim$ and each $f_j \in
\Hil$,
\begin{equation}
\label{AssUnif}
\|f_j\|_\infty \defn \sup_{x_j} |f_j(x_j)| \leq c \; \|f_j\|_\Hil.
\end{equation}
This condition is satisfied for many kernel classes including Sobolev
spaces, and any univariate kernel function\footnote{Indeed, we have
\begin{align*}
\sup_{x_j}|f_j(x_j)| & = \sup_{x_j}|\langle f_j(.), \Ker(.,x_j)
\rangle_{\Hil}| \; \leq \; \sup_{x_j} \sqrt{\Ker(x_j,x_j)}
\|f_j\|_\Hil.
\end{align*}
} bounded uniformly by $c$. Such a condition is routinely imposed for
proving upper bounds on rates of convergence for non-parametric least
squares in the univariate case $\pdim=1$ (see
e.g.~\cite{Sto85,vandeGeer}).  Note that this univariate boundedness
does not imply that the multivariate functions $f = \sum_{j \in S}
f_j$ in $\Fclass$ are uniformly bounded independently of $(\pdim,
\spindex)$; rather, they can take on values of the order $\sqrt{\s}$.

The following result applies to any class $\MyBigClass$ of sparse
additive models based on the univariate Hilbert space satisfying
condition~\eqref{AssUnif}, and to the estimator~\eqref{EqnMest} based
on $\numobs$ i.i.d. samples $(\exi, \yi)_{i=1}^{n}$ from the
observation model~\eqref{EqnLinObs}.
\btheos
\label{ThmPolyAchievable}
Let $\fhat$ be any minimizer of the convex program~\eqref{EqnMest}
with regularization parameters \mbox{$\regpar \geq 16 \gamcrit$} and
$\regpartwo \geq 16 \gamcrit^2$.  Then provided that \mbox{$\numobs
  \gamcrit^2 = \Omega(\log(1/\gamcrit))$,} there are universal
constants $(\UniCon, \plaincon_1, \plaincon_2)$ such that
\begin{align}
\label{EqnPolyAchievable}
\mprob \biggr[ \max \{ \|\fhat - \fstar\|_2^2, \; \|\fhat -
  \fstar\|_\numobs^2 \} \geq \UniCon \big \{ \spindex \regpar^2 +
  \spindex \regpartwo \big \} \biggr] & \leq \plaincon_1
\exp(-\plaincon_2 \numobs \gamcrit^2).
\end{align}

\etheos

\noindent We provide the proof of Theorem~\ref{ThmPolyAchievable} in
Section~\ref{SecAchResults}.

\paragraph{Remarks:}  

First, the technical condition $\numobs \gamcrit^2 =
\Omega(\log(1/\gamcrit))$ is quite mild, and satisfied in most cases
of interest, among them the kernels considered below in
Corollaries~\ref{CorAchieveFinite} and~\ref{CorAchieveSmooth}.

Second, note that setting $\regpar = c \gamcrit$ and $\regpartwo = c
\gamcrit^2$ for some constant $c \in (16, \infty)$ yields the rate
$\Theta(\s \gamcrit^2 + \s \regpartwo) = \Theta (\frac{\s \log d}{n} +
\s \unicrit^2)$.  This rate may be interpreted as the sum of a subset
selection term ($\frac{\spindex \log \pdim}{\numobs}$) and an
$\s$-dimensional estimation term ($\spindex \unicrit^2$). Note that
the subset selection term ($\frac{\spindex \log \pdim}{\numobs}$) is
independent of the choice of Hilbert space $\Hil$ whereas the
$\s$-dimensional estimation term is independent of the ambient
dimension $\pdim$.  Depending on the scaling of the triple $(\numobs,
\pdim, \s)$ and the smoothness of the univariate RKHS $\Hilb$, either
the subset selection term or function estimation term may dominate.
In general, if $\frac{\log \pdim}{n} = o(\unicrit^2)$, the
$\s$-dimensional estimation term dominates, and vice versa
otherwise. At the boundary, the scalings of the two terms are
equivalent.

Finally, for clarity, we have stated our result in the case where the
univariate Hilbert space $\Hil$ is identical across all co-ordinates.
However, our proof extends with only notational changes to the general
setting, in which each co-ordinate $j$ is endowed with a (possibly
distinct) Hilbert space $\Hil_j$.  In this case, the $M$-estimator
returns a function $\fhat$ such that (with high probability)
\begin{align*}
\max \big \{ \|\fhat - \fstar\|_\numobs^2, \; \|\fhat - \fstar\|_2^2
\big \} \; \; & \leq \; C \biggr \{ \frac{\spindex \log
  \pdim}{\numobs} + \sum_{j \in \Sset} \unicritpl^2_{\numobs, j}
\biggr \},
\end{align*}
where $\unicritpl_{\numobs, j}$ is the critical univariate rate
associated with the Hilbert space $\Hil_j$, and $\Sset$ is the subset
on which $\fstar$ is supported. \\
\vspace*{.1in}

Theorem~\ref{ThmPolyAchievable} has a number of corollaries, obtained
by specifying particular choices of kernels. First, we discuss
$\kerrank$-rank operators, meaning that the kernel function $\Ker$ can
be expanded in terms of $\kerrank$ eigenfunctions.  This class
includes linear functions, polynomial functions, as well as any
function class based on finite dictionary expansions. First we present
a corollary for finite-rank kernel classes.
\bcors
\label{CorAchieveFinite}
Under the same conditions as Theorem~\ref{ThmPolyAchievable}, consider
an univariate kernel with finite rank $\kerrank$.  Then any solution
$\fhat$ to the problem~\eqref{EqnMest} with $\regpar = c \gamcrit$ and
$\regpartwo = c \gamcrit^2$ with $16 \leq c < \infty$ satisfies
\begin{align}
\mprob \biggr[ \max \big \{ \|\fhat - \fstar\|_\numobs^2, \|\fhat -
  \fstar\|_2^2 \big \} \geq \UniCon \big \{ \frac{\spindex \log
    \pdim}{\numobs} + \spindex \frac{\kerrank}{\numobs} \big \}
  \biggr] & \leq \plaincon_1 \exp \big(-\plaincon_2 (\kerrank + \log
\pdim) \big).
\end{align}
\ecors
\spro
It suffices to show that the critical univariate
rate~\eqref{EqnCritCon} satisfies the scaling $\unicrit^2 =
\order(\kerrank/\numobs)$. For a finite-rank kernel and any $t > 0$,
we have
\begin{align*}
\LocRad(t) & = \frac{1}{\sqrt{\numobs}} \sqrt{\sum_{j=1}^\kerrank \min
  \{ t^2, \mu_j \} } \; \leq \; t \, \sqrt{\frac{\kerrank}{\numobs}},
\end{align*}
from which the claim follows by the definition~\eqref{EqnCritCon}.
\fpro

Next, we present a result for the RKHS's with infinitely many
eigenvalues, but whose eigenvalues decay at a rate $\mu_k \simeq
(1/k)^{2 \alpha}$ for some parameter $\alpha > 1/2$.  Among other
examples, this type of scaling covers the case of Sobolev spaces, say
consisting of functions with $\alpha$ derivatives
(e.g.,~\cite{BirSol67,Gu02}).
\bcors
\label{CorAchieveSmooth}
Under the same conditions as Theorem~\ref{ThmPolyAchievable}, consider
an univariate kernel with eigenvalue decay $\mu_k \simeq (1/k)^{2
  \alpha}$ for some $\alpha > 1/2$.  Then the kernel estimator defined
in \eqref{EqnMest} with $\regpar = c \gamcrit$ and $\regpartwo = c
\gamcrit^2$ with $16 \leq c < \infty$ satisfies
\begin{align}
\mprob \biggr[ \max \big\{ \|\fhat - \fstar\|_\numobs^2, \|\fhat -
  \fstar\|_2^2 \big \} \geq \UniCon \big \{ \frac{\spindex \log
    \pdim}{\numobs} + \spindex \big(\frac{1}{\numobs} \big)^{\frac{2
      \smooth}{2 \smooth +1}} \big \} \biggr] & \leq \plaincon_1 \exp
\big ( - \plaincon_2 (\numobs^{\frac{1}{2 \smooth +1}} + \log \pdim)
\big).
\end{align}
\spro
As in the previous corollary, we need to compute the critical
univariate rate $\unicrit$.  Given the assumption of polynomial
eigenvalue decay, a truncation argument shows that \mbox{$\LocRad(t) =
  \order \big(\frac{t^{1-\frac{1}{2 \smooth}}}{\sqrt{\numobs}}
  \big)$.}  Consequently, the critical univariate
rate~\eqref{EqnCritCon} satisfies the scaling $\unicrit^2 \asymp
\unicrit^{1-\frac{1}{2 \smooth}}/\sqrt{\numobs}$, or equivalently,
\mbox{$\unicrit^2 \asymp \numobs^{-\frac{2 \smooth}{2 \smooth +1}}$.}
\fpro
\ecors

\subsection{Minimax lower bounds}
\label{SecMiniLower}

In this section, we provide minimax lower bounds in $\LTP$ error so as
to complement the achievability results derived in
Theorem~\ref{ThmPolyAchievable}. Given the function class $\Fclass$,
the minimax $L^2(\mprob)$-error is given by
\begin{align}
\label{EqnDefnMinimax}
\MiniMaxPop & \defn \inf_{\fhat_\numobs} \sup_{\fstar \in \Fclass} \|
\fhat_\numobs - \fstar\|_2^2,
\end{align}
where the infimum is taken over all measurable functions of $\numobs$
samples $\{(\exi, \yi)\}_{i=1}^\numobs$.  As defined, this minimax
error is a random variable, and our goal is to obtain a lower bound in
probability.

Central to our proof of the lower bounds is the metric entropy
structure of the univariate reproducing kernel Hilbert spaces. More
precisely, our lower bounds depend on the \emph{packing entropy,}
defined as follows.  Let $(\GenSet, \rho)$ be a totally bounded metric
space, consisting of a set $\GenSet$ and a metric $\rho: \GenSet
\times \GenSet \rightarrow \mathbb{R}_+$.  An $\epsilon$-packing of
$\GenSet$ is a collection $\{\f^1, \ldots, \f^\PackNum \} \subset
\GenSet$ such that $\rho(\f^i, \f^j) \geq \epsilon$ for all $i \neq
j$.  The $\epsilon$-packing number $\PackNum(\epsilon; \GenSet, \rho)$
is the cardinality of the largest $\epsilon$-packing.  The packing
entropy is the simply the logarithm of the packing number, namely the
quantity $\PackEnt(\epsilon; \GenSet, \rho)$, to which we also refer
as the metric entropy.  In this paper, we derive explicit minimax
lower bounds for two different scalings of the univariate metric
entropy. \\

\paragraph{Logarithmic metric entropy:} There exists
some $\mdim > 0$ such that
\begin{align}
\label{EqnLogMet}
\log M(\epsilon; \Ball_\Hil(1), L^2(\mathbb{P})) & \simeq \mdim \,
\log(1/\epsilon) \qquad \mbox{for all $\epsilon \in (0,1)$.}
\end{align}
Function classes with metric entropy of this type include linear
functions (for which $\mdim = k$), univariate polynomials of degree
$k$ (for which $\mdim = k + 1$), and more generally, any function
space with finite VC-dimension~\cite{vanderVaart96}.  This type of
scaling also holds for any RKHS based on a kernel with rank $\mdim$
(e.g., see~\cite{CarlTrieb80}), and these finite-rank kernels include
both linear and polynomial functions as special cases.

\paragraph{Polynomial metric entropy} There exists some
$\smooth > 0$ such that
\begin{align}
\label{EqnPolyMet}
\log M(\epsilon; \Ball_\Hil(1), L^2(\mathbb{P})) & \simeq
(1/\epsilon)^{1/\smooth} \qquad \mbox{for all $\epsilon \in (0,1)$.}
\end{align}
Various types of Sobolev/Besov classes exhibit this type of metric
entropy decay~\cite{BirSol67,Gu02}.  In fact, any RKHS in which the
kernel eigenvalues decay at a rate $k^{-2\NonPar}$ have a metric
entropy with this scaling~\cite{Carl90,CarlTrieb80}. \\

\noindent We are now equipped to state our lower bounds on the minimax
risk~\eqref{EqnDefnMinimax}:
\btheos
\label{ThmLower}
Given $\numobs$ i.i.d. samples from the sparse additive
model~\eqref{EqnLinObs} with sparsity $\s \leq \pdim/4$, there is an
universal constant $\UniCon > 0$ such that:
\begin{enumerate}
\item[(a)] For a univariate class $\Hilb$ with logarithmic metric
  entropy~\eqref{EqnLogMet} indexed by parameter $\ParDim$, we have
\begin{align}
\label{EqnLogLower}
\MiniMaxPop & \; \geq \; \UniCon \biggr \{ \frac{\spindex \log
  (\pdim/\spindex)}{\numobs} \, + \, \spindex \,
\frac{\ParDim}{\numobs } \biggr \}
\end{align}
with probability greater than $1/2$.
\item[(b)] For a univariate class $\Hilb$ with polynomial metric
  entropy~\eqref{EqnPolyMet} indexed by $\smooth$, we have
\begin{align}
\label{EqnPolyLower}
\MiniMaxPop & \; \geq \; \UniCon \, \biggr \{ \frac{\spindex \log
  (\pdim/\spindex)}{\numobs} \, + \, \spindex \,\big(\frac{1}{\numobs}
\big)^{\frac{2 \smooth}{2 \smooth +1}} \biggr \}
\end{align}
with probability greater than $1/2$.
\end{enumerate}
\etheos

\noindent The proof of Theorem~\ref{ThmLower} is provided in
Section~\ref{SecLowerBound}.  Our choice of stating bounds that hold
with probability $1/2$ is simply a convention often used in
information-theoretic approaches (see, for instance, the
papers~\cite{Hasminskii78,YanBar99,Yu}). We note that analogous lower
bounds can established with probabilities arbitrarily close to one,
albeit at the expense of worse constants.  The most important
consequence of Theorem~\ref{ThmLower} is in establishing the
minimax-optimality of the results given in
Corollary~\ref{CorAchieveFinite} and~\ref{CorAchieveSmooth}; in
particular, in the regime of sub-linear sparsity (i.e., for which
$\log \pdim = \order(\log (\pdim/\s))$), the combination of
Theorem~\ref{ThmLower} with these corollaries identifies the minimax
rates up to constant factors.

\subsection{Comparison with other estimators}
\label{SecCompare}

It is interesting to compare these convergence rates in $\LTPn$ error
with those established in past
work~\cite{KolYua08,Meier09,RavLiuLafWas08} using different
estimators. Ravikumar et. al~\cite{RavLiuLafWas08} show that any
solution to their back-fitting method is consistent in terms of
mean-squared error risk (see Theorem 3 in their paper) but they don't
have rates in the regime where $\s \rightarrow \infty$.  An earlier
method of Koltchinskii and Yuan~\cite{KolYua08} is based regularizing
the least-squares loss with the $(\Hilb, 1)$-norm penalty---that is,
$\sum_{j=1}^{\pdim}{\|f_j\|_{\Hilb}}$ but no $(\|.\|_n, 1)$-norm
penalty; Theorem 2 in their paper presents a rate that captures the
decomposition into two terms, a subset selection and $\s$-dimensional
estimation term. In quantitative terms however, their rates are looser
than those given here; in particular, their bound includes a term of
the order $\frac{\s^3 \log \pdim}{n}$, which is larger than the bound
in Theorem~\ref{ThmPolyAchievable}.  Meier et al.~\cite{Meier09}
analyze a different $M$-estimator to the one we analyze in this
paper. Rather than adding two separate $(\Hilb, 1)$-norm and an
$(\|.\|_n, 1)$-norm penalties, they combine the two terms into a
single sparsity and smoothness penalty. For their estimator, Meier et
al.~\cite{Meier09} establish a convergence rate of the form $\order(\s
(\frac{\log \pdim}{n})^{\frac{2 \NonPar}{2 \NonPar+ 1}}\big)$ in the
case of $\NonPar$-smooth Sobolev spaces (see Theorem 1 in their
paper). This result is sub-optimal compared to the optimal rate proven
in Theorem~\ref{ThmLower}(b). More precisely, we either have
$\frac{\log \pdim}{\numobs} < (\frac{\log \pdim}{\numobs})^{\frac{2
    \NonPar}{2 \NonPar+ 1}}$, when subset selection term dominates, or
\mbox{$(\frac{1}{n})^{\frac{2 \NonPar}{2 \NonPar+ 1}} < (\frac{\log
    \pdim}{\numobs})^{\frac{2 \NonPar}{2 \NonPar+ 1}}$,} when the
$\s$-dimensional estimation term dominates.  In all of the
above-mentioned methods, it is unclear whether or not sharper analysis
would yield better rates. Koltchinskii and Yuan~\cite{KolYua10Journal}
analyzes the same estimator as the $M$-estimator~\eqref{EqnMest} and
achieve the same rates as in Theorem~\ref{ThmPolyAchievable}, under a
global boundedness condition. In the following section, we show that
rates in their paper are not minimax optimal for Sobolev spaces when
$\s = \Omega(\sqrt{n})$.

\subsection{Upper bounds under a global boundedness assumption}
\label{SecKYCompare}

As discussed previously in the introduction, the past work of
Koltchinski and Yuan~\cite{KolYua10Journal}, referred to as KY for
short, is based on the $M$-estimator~\eqref{EqnMest}. In terms of
rates obtained, they establish a convergence rate based on two terms
as in Theorem~\ref{ThmPolyAchievable}, but with a pre-factor that
depends on the global quantity
\begin{align}
\label{EqnGlobalBound}
\GLOB = \sup_{f \in \MyBigClass} \|f\|_\infty \; = \; \sup_{f \in
  \MyBigClass} \sup_x \; |f(x)|,
\end{align}
assumed to be bounded independently of dimension and sparsity.  Such
types of global boundedness conditions are fairly standard in
classical non-parametric estimation and is equivalent to univariate
boundedness, and they have no effect on minimax rates.  In sharp
contrast, the analysis of this section shows that for sparse additive
models in the regime $\s = \Omega(\sqrt{n})$, such global boundedness
can \emph{substantially speed up} minimax rates, showing that the
rates proven in KY are not minimax optimal for these classes. The
underlying insight is as follows: when the sparsity grows, imposing
global boundedness over $\s$-variate functions substantially reduces
the effective dimension from its original size $\s$ to a lower
dimensional quantity, which we denote by $\s K_\GLOB(\s,n)$, and
moreover, the quantity $K_\GLOB(\s,n) \rightarrow 0$ when $\s =
\Omega(\sqrt{n})$ as described below.

Recall the definition~\eqref{EqnDefnFclass} of the function class
$\MyBigClass$.  The model considered in the KY paper is the smaller
function class
\begin{align*}
\MyNewClass & \defn \bigcup_{\substack{\Sset \subset \{1, 2, \ldots,
    \pdim\} \\ |\Sset| = \kdim}} \MyNewSmallClass,
\end{align*}
where $\MyNewSmallClass \defn \big \{ f = \sum_{j \in \Sset} f_j \,
\mid \, f_j \in \Hil, \mbox{ and } f_j \in \Ball_{\Hil}(1) \; \,
\forall \; j \in \Sset\; \mbox{and}\;\|f\|_{\infty} \leq B \big \}$.

The following theorem provides sharper rates for the Sobolev case, in
which each univariate Hilbert space has eigenvalues decaying as $\mu_k
\simeq k^{-2 \alpha}$ for some smoothness parameter $\alpha > 1/2$.
Our probabilistic bounds involve the quantity
\begin{equation}
\label{EqnDefnSmallRate}
 \UnifBoundRate \defn \max \big(\sqrt{\frac{\s \log (\pdim/\s)}{n}},
 B^{1/2}(\frac{\s^{\frac{1}{\alpha}}\log \s}{n})^{1/4} \big),
\end{equation}
and our rates are stated in terms of the function
\begin{align}
\label{EqnBfun}
K_\GLOB(\s, \numobs) & \defn \GLOB \sqrt{\log \s} (\s^{-1/2\alpha}
\numobs^{1/(4 \alpha + 2)})^{2 \alpha - 1}.
\end{align}
Note that $K_\GLOB(\s, \numobs) \rightarrow 0$ if $\s = \Omega(\sqrt{n})$. With this notation, we have the following \emph{upper bound} on the minimax risk over the function class $\MyNewClass$.
\btheos
\label{ThmFast}
Consider a Sobolev RKHS $\Hil$ with eigenvalue decay $k^{-2 \alpha}$ and
eigenfunctions such that $\|\phi_k\|_{\infty} \leq C < \infty$.  Then
there are universal constants $(c_1,c_2, \kappa)$ such that
with probability greater than $1 - 2 \exp \big( - c_1 n \UnifBoundRate^2
\big)$, we have
\begin{align}
%
\label{EqnFastPop}
\underbrace{\min_{\hat{f}} \max_{f^* \in \MyNewClass}\|\fhat -
  \fstar\|_2^2}_{\MiniMaxFast} & \leq \kappa^2 (1+B) C \s n^{-\frac{2
    \alpha}{2 \alpha + 1}}\biggr(K_B(s,n)+ n^{-1/(2 \alpha + 1)} \log
(\pdim/\s) \biggr),
\end{align}
as long as $\numobs \UnifBoundRate^2 =
\Omega(\log(1/\UnifBoundRate))$.  
\etheos

We provide the proof of Theorem~\ref{ThmFast} in
Section~\ref{SecProofFast}; it is based on analyzing directly the
least-squares estimator over $\MyNewClass$.  The assumption that
$\|\phi_k\|_{\infty} \leq C < \infty$ for all $k$ includes the usual
Sobolev spaces in which $\phi_k$ are (rescaled) Fourier basis
functions.  An immediate consequence of Theorem~\ref{ThmFast} is that
the minimax rates over the function class $\MyNewClass$ can be
strictly faster than minimax rates for the class $\MyBigClass$ (which
does not assume global boundedness).  Recall that the minimax lower
bound from Theorem~\ref{ThmLower} (b) takes the form:
\begin{equation*}
\MiniMaxPop \; \defn \; \min_{\fhat} \max_{f^* \in \MyBigClass}\|\fhat
- \fstar\|_2^2 \geq \UniCon_1 \s n^{-\frac{2 \alpha}{2 \alpha +
    1}}\biggr(1+ n^{-1/(2 \alpha + 1)} \log (\pdim/\s) \biggr),
\end{equation*}
for a universal constant $\UniCon_1$. Note that up to constant
factors, the achievable rate~\eqref{EqnFastPop} from
Theorem~\ref{ThmFast} is the same except that the term $1$ is replaced
by the function $K_\GLOB(\s, \numobs)$.  Consequently, for scalings of
$(\s, \numobs)$ such that $K_\GLOB(\s, \numobs) \rightarrow 0$, global
boundedness conditions lead to strictly faster rates.

\bcors
\label{CorFast}
Under the conditions of Theorem~\ref{ThmFast}, we have
\begin{align*}
\frac{\MiniMaxPop}{\MiniMaxFast} & \geq \frac { \UniCon_1 (1+ n^{-1/(2
    \alpha + 1)} \log (\pdim/\s))} {C \, \kappa^2(1+ \GLOB) \,
  (K_B(s,n)+ n^{-1/(2 \alpha + 1)} \log (\pdim/\s))} \; \rightarrow
+\infty
\end{align*}
with probability at least $1/2$, whenever $B = \order(1)$ and
$K_\GLOB(\s, \numobs) \rightarrow 0$.
\ecors

\paragraph{Remarks:}
The quantity $K_\GLOB(\s, \numobs)$ is guaranteed to decay to zero as
long as the sparsity index $\s$ grows in a non-trivial way with the
sample size.  For instance, if we have $\kdim = \Omega(\sqrt{n})$ for
a problem of dimension $\pdim = \order(n^{\beta})$ for any $\beta \geq
1/2$, then it can be verified that $f_\GLOB(\s, \numobs) = o(1)$.  As
an alternative view of the differences, it can be noted that there are
scalings of $(\numobs, \s, \pdim)$ for which the minimax rate
$\MiniMaxPop$ over $\MyBigClass$ is constant---that is, does not
vanish as $\numobs \rightarrow +\infty$---while the minimax rate
$\MiniMaxFast$ does vanish.  As an example, consider the Sobolev class
with smoothness $\alpha = 2$, corresponding to twice-differentiable
functions.  For a sparsity index $\kdim = \order(\numobs^{4/5})$, then
Theorem~\ref{ThmLower}(b) implies that $\MiniMaxPop = \Omega(1)$, so
that the minimax rate over $\MyBigClass$ is strictly bounded away from
zero for all sample sizes.  In contrast, under a global boundedness
condition, Theorem~\ref{ThmFast} shows that the minimax rate is upper
bounded as $\MiniMaxFast = \order \big(n^{-1/5} \sqrt{\log n}\big)$,
which tends to zero. \\

In summary, Theorem~\ref{ThmFast} and Theorem~\ref{ThmLower}(b)
together show that the minimax rates over $\MyBigClass$ and
$\MyNewClass$ can be drastically different.  Thus, global boundedness
is a stringent condition in the high-dimensional setting; in
particular, the rates given in Theorem 3 of Koltchinskii and
Yuan~\cite{KolYua10Journal} are not minimax optimal when $s =
\Omega(\sqrt{\numobs})$. 


\section{Proofs}
\label{SecProofs}

In this section, we provide the proofs of our three main theorems. For
clarity in presentation, we split the proofs up into a series of
lemmas, with the bulk of the more technical proofs deferred to the
appendices.  This splitting allows our presentation in
Section~\ref{SecProofs} to be relatively streamlined.


\subsection{Proof of Theorem~\ref{ThmPolyAchievable}}
\label{SecAchResults}

At a high-level, Theorem~\ref{ThmPolyAchievable} is based on an
appropriate adaptation to the non-parametric setting of various
techniques that have been developed for bounding the error to those in
sparse linear regression (e.g.,~\cite{BicRitTsy08,Neg09}).  In
contrast to the parametric setting where classical tail bounds are
sufficient, controlling the error terms in the non-parametric case
requires more advanced techniques from empirical process theory. In
particular, we make use of concentration theorems for Gaussian and
empirical processes
(e.g.,~\cite{Chung06,Ledoux01,Massart00b,Pisier89,vandeGeer}) as well
as results on the Rademacher complexity of kernel
classes~\cite{Bartlett05,Mendelson02}.

At the core of the proof are three technical lemmas.  First,
Lemma~\ref{LemNetherlands} provides an upper bound on the Gaussian
complexity of any function of the form $f = \sum_{j=1}^\pdim f_j$ in
terms of the norms $\|\cdot\|_{\Hil,1}$ and $\|\cdot\|_{\numobs,1}$
previously defined.  Lemma~\ref{LemConeBound} exploits the notion of
decomposability~\cite{Neg09}, as applied to these norms, in order to
show that the error function belongs to a particular cone-shaped set.
Finally, Lemma~\ref{LemLowerBound} establishes an upper bound on the
$\LTP$ error of our estimator in terms of the $\LTPn$ error which
completes the proof. The latter lemma can be interpreted as proving
that our problem satisfies non-parametric analog of a restricted
strong convexity~\cite{Neg09} or restricted eigenvalue
condition~\cite{BicRitTsy08}. The proof of Lemma~\ref{LemLowerBound}
involves a new approach that combines the Sudakov
minoration~\cite{Pisier89} with a one-sided concentration bound for
non-negative random variables~\cite{Chung06}.

Throughout the proof, we use $\UniCon$ and $\plaincon_i$, $i=1,2,3,4$
to denote universal constants, independent of $(\numobs, \pdim,
\spindex)$. Note that the precise numerical values of these constants
may change from line to line. The reader should recall the definitions
of $\unicrit$ and $\gamcrit$ from equations~\eqref{EqnCritCon}
and~\eqref{EqnGamCrit} respectively. For a subset \mbox{$A \subseteq
  \{ 1,2, \ldots, \pdim\}$} and a function of the form $f =
\sum_{j=1}^{\pdim} {f_j}$, we adopt the convenient notation
\begin{equation}
\label{EqnDecompNotation}
\Lpeone{f_A}\; \defn \; \sum_{j \in A} \|f_j\|_\numobs, \quad
\mbox{and} \quad \Hilone{f_A} \; \defn \; \sum_{j \in A} \|f_j\|_\Hil.
\end{equation}

We begin by establishing an inequality on the error function
$\DeltaHat \defn \fhat - \fstar$.  Since $\fhat$ and $\fstar$ are,
respectively, optimal and feasible for the problem~\eqref{EqnMest}, we
are guaranteed that $\CostFun(\fhat) \leq \CostFun(\fstar)$, and hence
that the error function $\DeltaHat$ satisfies the bound
\begin{equation*}
\frac{1}{2n}\sum_{i=1}^\numobs {(\wi + \mu - \ybar -
  \DeltaHat(\exi))^2}+ \regpar \Lpeone{\fhat} + \regpartwo
\Hilone{\fhat} \leq \frac{1}{2n}\sum_{i=1}^\numobs {(\wi + \mu -
  \ybar)^2}+ \regpar \Lpeone{\fstar} + \regpartwo \Hilone{\fstar}.
\end{equation*}
Some simple algebra yields the bound
\begin{equation}
\label{EqRiskBound}
\frac{1}{2}\|\DeltaHat\|_\numobs^2 \leq \big | \frac{1}{\numobs}
\sum_{i=1}^\numobs {\wi\DeltaHat(\exi)} \big | + |\ybar-\mu | \big |
\frac{1}{\numobs} \sum_{i=1}^\numobs {\DeltaHat(\exi)} \big |+ \regpar
\Lpeone{\DeltaHat} + \regpartwo \Hilone{\DeltaHat}.
\end{equation}
Following the terminology of van de Geer~\cite{vandeGeer}, we refer to
this bound as our \emph{basic inequality}.

\subsubsection{Controlling deviation from the mean}

Our next step is to control the error due to estimating the mean
$|\ybar-\mu|$.  We begin by observing that this error term can be
written as $\ybar - \mu = \frac{1}{\numobs} \sum_{i=1}^{\numobs} (y_i
- \mu)$, Now consider the random variable $y_i - \mu = \sum_{j \in
  S}{f_j^*(\xij)} + \wi$. It is the sum of the $\s$ independent random
variables $f_j^*(\xij)$, each bounded in absolute value with one,
along with the independent sub-Gaussian noise term $\wi$;
consequently, the variable $y_i - \mu$ is sub-Gaussian with parameter
at most $\sqrt{ \s + 1}$ (see, e.g., Lemma 1.4 in Buldygin and
Kozachenko~\cite{BulKoz}).  By applying standard sub-Gaussian tail
bounds, we have $\mP (|\ybar-\mu| > t) \leq 2 \exp(- \frac{\numobs
  t^2}{2 (\s+1)})$, and hence, if we define the event
$\mathcal{C}(\gamma_n) = \{ |\ybar-\mu| \leq \sqrt{\s}\gamma_n \}$, we
are guaranteed
\begin{equation*}
\mP[\mathcal{C}(\gamma_n)] \geq 1 - 2 \exp(- \frac{n \gamma_n^2}{4}).	
\end{equation*}	
For the remainder of the proof, we condition on the event
$\mathcal{C}(\gamma_n)$.  Under this conditioning, the
bound~\eqref{EqRiskBound} simplifies to:
\begin{equation*}
\frac{1}{2}\|\DeltaHat\|_\numobs^2 \leq \big | \frac{1}{\numobs}
\sum_{i=1}^\numobs {\wi\DeltaHat(\exi)} \big | + \sqrt{\s} \gamcrit
\|\DeltaHat\|_n + \regpar \Lpeone{\DeltaHat} + \regpartwo
\Hilone{\DeltaHat},
\end{equation*}
where we have applied the Cauchy-Schwarz inequality to write $\big |
\frac{1}{\numobs} \sum_{i=1}^\numobs {\DeltaHat(\exi)} \big | \leq
\|\DeltaHat\|_\numobs$.

\subsubsection{Controlling the Gaussian complexity term}

The following lemma provides control the Gaussian complexity term on
the right-hand side of inequality~\eqref{EqRiskBound} by bounding the
Gaussian complexity for the univariate functions $\DeltaHat_j$, $j =
1, 2, \ldots, \pdim$ in terms of their $\| \cdot \|_\numobs$ and $\|
\cdot \|_\Hil$ norms. In particular, recalling that $\gamcrit = \kappa
\max \{ \sqrt{\frac{\log \pdim}{\numobs}}, \; \unicrit \}$, we have
the following lemma.

\blems
\label{LemNetherlands}
Define the event
\begin{align}
\TailEvent(\gamcrit) & \defn \biggr \{ \forall \;\; j=1, 2, \ldots,
\pdim,\; \big| \frac{1}{\numobs} \sum_{i=1}^\numobs \wi
\DeltaHat_j(\xij) \big| \leq 8 \gamcrit^2 \;
\|\DeltaHat_j\|_\Hil + 8 \gamcrit \; \|\DeltaHat_j\|_\numobs
\biggr \}.
\end{align}
Then under the condition $\numobs \gamcrit^2 =
\Omega(\log(1/\gamcrit))$, we have
\begin{align}
\label{EqnNetherlands}
\mathbb{P}(\TailEvent(\gamcrit)) \geq 1 - \plaincon_1
\exp(-\plaincon_2 \numobs \gamcrit^2).
\end{align}
\elems
\noindent
The proof of this lemma, provided in Appendix~\ref{AppLemNetherlands},
uses concentration of measure for Lipschitz functions over Gaussian
random variables~\cite{Ledoux01} combined with peeling and weighting
arguments~\cite{Alex87,vandeGeer}. In particular, the subset selection
term $(\frac{\s \log \pdim}{\numobs})$ in
Theorem~\ref{ThmPolyAchievable} arises from taking the maximum over
all $\pdim$ components.

The remainder of our analysis involves conditioning on the event
$\TailEvent(\gamcrit)\cap \mathcal{C}(\gamma_n)$. Using
Lemma~\ref{LemNetherlands}, when conditioned on the event
$\TailEvent(\gamcrit) \cap \mathcal{C}(\gamma_n)$ we have:
\begin{eqnarray}
\label{EqnDecomp} 
\|\DeltaHat\|_\numobs^2 & \leq & 2 \sqrt{\s} \gamcrit
\|\DeltaHat\|_\numobs + (16 \gamcrit + 2 \regpar) \Lpeone{\DeltaHat} +
(16 \gamcrit^2 + 2 \regpartwo) \Hilone{\DeltaHat}.
\end{eqnarray}

\subsubsection{Exploiting decomposability}

Recall that $S$ denotes the true support of the unknown function
$\fstar$.  By the definition~\eqref{EqnDecompNotation}, we can write
$\Lpeone{\DeltaHat} = \Lpeone{\DeltaHat_S} +
\Lpeone{\DeltaHat_{S^c}}$, where $\DeltaHat_S \defn \sum_{j \in
  S}{\DeltaHat_j}$ and $\DeltaHat_{S^c} \defn \sum_{j \in
  S^c}{\DeltaHat_j}$.  Similarly, we have an analogous representation
for $\Hilone{\DeltaHat}$.  The next lemma shows that conditioned on
the event $\TailEvent(\gamcrit)$, the quantities
$\|\DeltaHat\|_{\Hil,1}$ and $\|\DeltaHat\|_{\numobs,1}$ are not
significantly larger than the corresponding norms as applied to the
function $\DeltaHat_\Sset$.

\blems
\label{LemConeBound}
Conditioned on the events $\TailEvent(\gamcrit)$ and
$\mathcal{C}(\gamma_n)$, and with the choices $\regpar \geq 16
\gamcrit$ and $\regpartwo \geq 16 \gamcrit^2$, we have
\begin{align}
\label{EqnConeBound}
\regpar \Lpeone{\DeltaHat} + \regpartwo \Hilone{\DeltaHat}
& \leq 4 \regpar \Lpeone{\DeltaHat_\Sset} + 4 \regpartwo
\Hilone{\DeltaHat_\Sset} + \frac{1}{2}\s \gamcrit^2.
\end{align}
\elems
\noindent The proof of this lemma, provided in
Appendix~\ref{AppLemConeBound}, is based on the
decomposability~\cite{Neg09} of the $\|\cdot\|_{\Hil, 1}$ and
$\|\cdot\|_{\numobs, 1}$ norms. This lemma allows us to exploit the
sparsity assumption, since in conjunction with
Lemma~\ref{LemNetherlands}, we have now bounded the right-hand side of
the inequality~\eqref{EqnDecomp} by terms involving only
$\DeltaHat_\Sset$.  

For the remainder of the proof of Theorem~\ref{ThmPolyAchievable}, we
assume $\regpar \geq 16 \gamcrit$ and $\regpartwo \geq 16
\gamcrit^2$. In particular, still conditioning on
$\mathcal{C}(\gamma_n) \cap \TailEvent(\gamcrit)$ and applying
Lemma~\ref{LemConeBound} to inequality~\eqref{EqnDecomp}, we obtain
\begin{eqnarray*}
\|\DeltaHat\|_\numobs^2 & \leq & 2 \sqrt{\s} \gamcrit
\|\DeltaHat\|_\numobs + 3 \regpar \Lpeone{\DeltaHat} + 3 \regpartwo
\Hilone{\DeltaHat}\\
& \leq & 2 \sqrt{\s}\regpar \|\DeltaHat\|_\numobs + 12 \regpar
\Lpeone{\DeltaHat_\Sset} + 12 \regpartwo \Hilone{\DeltaHat_\Sset} +
\frac{3}{32}\s \regpartwo,
\end{eqnarray*}
Finally, since both $\fhat_j$ and $\fstar_j$ belong to
$\Ball_\Hil(1)$, we have $\|\DeltaHat_j\|_\Hil \leq \|\fhat_j\|_\Hil +
\|\fstar_j\|_\Hil \leq 2$, which implies that
$\Hilone{\DeltaHat_\Sset} \leq 2 \spindex$, and hence
\begin{align}
\label{EqnUpperBound}
\|\DeltaHat\|_\numobs^2 & \leq 2 \sqrt{\s}\regpar
\|\DeltaHat\|_\numobs + 12 \regpar \Lpeone{\DeltaHat_\Sset} + 25 \s
\regpartwo.
\end{align}

\subsubsection{Upper bounding $\Lpeone{\DeltaHat_\Sset}$}

The final step is to control the term $\|\DeltaHat_\Sset\|_{\numobs,1}
= \sum_{j \in \Sset} \|\DeltaHat_j\|_\numobs$ that appears in the
upper bound~\eqref{EqnUpperBound}. Ideally, we would like to upper
bound it by a quantity of the order $\sqrt{\spindex}
\|\DeltaHat_S\|_2$. Such an upper bound would follow immediately if it
were phrased in terms of the population $\|\cdot\|_2$-norm rather than
the empirical-$\|\cdot\|_\numobs$ norm, but there are additional
cross-terms with the empirical norm.  Accordingly, a somewhat more
delicate argument is required, which we provide here.  First define
the events
\begin{align*}
\AuxEvent_j(\regpar) \defn \{\|\DeltaHat_j\|_\numobs & \leq 2
\|\DeltaHat_j\|_2 + \regpar\},
\end{align*}
and $\AuxEvent(\regpar) = \cap_{j=1}^\pdim \AuxEvent_j(\regpar)$. By
applying Lemma~\ref{LemGeneral} from Appendix~\ref{AppLemEquivNorms}
with $t = \regpar \geq 16 \gamcrit$ and $b=2$, we conclude that
$\|f_j\|_\numobs \leq 2 \|f_j\|_2 + \regpar$.  with probability
greater than $1 - \plaincon_1 \exp(-\plaincon_2 \numobs \regpar^2)$.
Consequently, if we define the event $\AuxEvent(\regpar) = \cap_{j \in
  S} \AuxEvent_j(\regpar)$, then this tail bound together with the
union bound implies that
\begin{align}
\label{EqnAuxBound}
\mprob[\AuxEvent^c(\regpar)] & \leq \s \; \plaincon_1
\exp(-\plaincon_2 \numobs \regpar^2) \; \leq \; \plaincon_1
\exp(-\plaincon'_2 \numobs \regpar^2),
\end{align}
where we have used the fact that $\regpar = \Omega(\sqrt{\frac{\log
    \s}{\numobs}})$.  Now, conditioned on the event
$\AuxEvent(\regpar)$, we have
\begin{eqnarray}
\label{EqnSquirrel}
\Lpeone{\DeltaHat_\Sset} & \leq & \sum_{j \in \Sset}
\|\DeltaHat_j\|_\numobs \; \leq \; 2 \sum_{j \in \Sset}
\|\DeltaHat_j\|_{2} + \s \regpar \; \leq \; 2 \sqrt{\s}
\|\DeltaHat_\Sset\|_2 + \s \regpar \; \leq \; 2 \sqrt{\s}
\|\DeltaHat\|_2 + \s \regpar.
\end{eqnarray}
Substituting this upper bound~\eqref{EqnSquirrel} on
$\Lpeone{\DeltaHat_\Sset}$ into our earlier
inequality~\eqref{EqnUpperBound} yields
\begin{eqnarray}
\label{EqnBasicTwo}
\|\DeltaHat\|_\numobs^2 & \leq & 2 \sqrt{\s}\regpar
\|\DeltaHat\|_\numobs + 24 \sqrt{\s} \regpar \|\DeltaHat\|_2 + 12 \s
\regpar^2 + 25 \s \regpartwo.
\end{eqnarray}

At this point, we encounter a challenge due to the unbounded nature of
our function class.  In particular, if $\|\DeltaHat\|_2$ were upper
bounded by $C\max(\|\DeltaHat\|_n, \sqrt{\s}\regpar, \sqrt{\s
  \regpartwo})$, then the upper bound~\eqref{EqnBasicTwo} would
immediately imply the claim of Theorem~\ref{ThmPolyAchievable}.  If
one were to assume global boundedness of the multivariate functions
$\fHat$ and $f^*$, as done in past work~\cite{KolYua10Journal}, then
an upper bound on $\|\DeltaHat\|_2$ of this form would directly follow
from known results (e.g., Theorem 2.1 in Bartlett et
al.~\cite{Bartlett05}.) However, since we do not impose global
boundedness, we need to develop a novel approach to obtaining suitable
control $\|\DeltaHat\|_2$, the task to which we now turn.

\subsubsection{Controlling $\|\DeltaHat\|_2$ for unbounded classes}

For the remainder of the proof, we condition on the event
$\AuxEvent(\regpar)\cap\TailEvent(\gamcrit)\cap
\mathcal{C}(\gamma_n)$.  We split our analysis into three
cases. Throughout the proof, we make use of the quantity
\begin{align}
\label{EqnDefnTildecrit}
\tildelcrit & \defn B \max(\sqrt{\s}\regpar, \sqrt{\s \regpartwo}),
\end{align}
where $B \in (1, \infty)$ is a constant to be chosen later in the
argument.
 
\paragraph{Case 1:} If $\|\DeltaHat\|_2 < \|\DeltaHat\|_n$, 
then combined with inequality~\eqref{EqnBasicTwo}, we conclude that
\begin{align*}
\|\DeltaHat\|_\numobs^2 \leq 2 \sqrt{\s}\regpar \|\DeltaHat\|_\numobs
+ 24 \sqrt{\s} \regpar \|\DeltaHat\|_n + 12 \s \regpar^2 + 25 \s
\regpartwo.
\end{align*}
This is a quadratic inequality in terms of the quantity
$\|\DeltaHat\|_\numobs$, and some algebra shows that it implies the
bound $\|\DeltaHat\|_\numobs \leq 15 \max(\sqrt{\s}\regpar, \sqrt{\s
  \regpartwo})$.  By assumption, we then have $\|\DeltaHat\|_2 \leq 15
\max(\sqrt{\s}\regpar, \sqrt{\s \regpartwo})$ as well, thereby
completing the proof of Theorem~\ref{ThmPolyAchievable}.


\paragraph{Case 2:} If $\|\DeltaHat\|_2 < \tildelcrit$, then together
with the bound~\eqref{EqnBasicTwo}, we conclude that
\begin{align*}
\|\DeltaHat\|_\numobs^2 \leq 2 \sqrt{\s}\regpar \|\DeltaHat\|_\numobs
+ 24 \sqrt{\s} \regpar \tilde{\delta_n} + 12 \s \regpar^2 + 25 \s
\regpartwo.
\end{align*}
This inequality is again a quadratic in $\|\DeltaHat\|_\numobs$;
moreover, note that by definition~\eqref{EqnDefnTildecrit} of
$\tildelcrit$, we have $\s \regpar ^2 + \s \regpartwo =
\order(\tildelcrit^2)$.  Consequently, this inequality implies that
$\|\DeltaHat\|_\numobs \leq C \tildelcrit$ for some constant $C$.  Our
starting assumption implies that $\|\DeltaHat\|_2 \leq \tildelcrit$,
so that the claim of Theorem~\ref{ThmPolyAchievable} follows in this
case.


\paragraph{Case 3:} Otherwise, we may assume that 
$\|\DeltaHat\|_2 \geq \tildelcrit$ and $\|\DeltaHat\|_2 \geq
\|\DeltaHat\|_\numobs$.  In this case, the
inequality~\eqref{EqnBasicTwo} together with the bound
$\|\DeltaHat\|_2 \geq \|\DeltaHat\|_\numobs$ implies that
\begin{align}
\label{EqnJunco}
\|\DeltaHat\|_\numobs^2 \leq 2 \sqrt{\s}\regpar \|\DeltaHat\|_2 + 24
\sqrt{\s} \regpar \|\DeltaHat\|_2 + 12 \s \regpar^2 + 25 \s
\regpartwo.
\end{align}
Our goal is to establish a lower bound on the left-hand-side---namely,
the quantity $\|\DeltaHat\|_\numobs^2$---in terms of
$\|\DeltaHat\|_2^2$.  In order to do so, we consider the function
class $\Gclass(\regpar, \regpartwo)$ defined by functions
of the form $g = \sum_{j=1}^\pdim g_j$, and such that
\begin{subequations}
\label{EqnSuperGclass}
\begin{align}
\label{EqnGclassCone}
\regpar \Lpeone{g} + \regpartwo \Hilone{g} & \leq \; 4 \regpar
\Lpeone{g_\Sset} + 4 \regpartwo \Hilone{g_\Sset} + \frac{1}{32}\s
\regpartwo, \\
\label{EqnGclassSet}
\|g_\Sset\|_{1, \numobs} & \leq \; 2 \sqrt{\s} \|g_\Sset\|_2 + \s
\regpar \qquad \mbox{and} \\
\|g\|_\numobs & \leq \|g\|_2.
\end{align}
\end{subequations}
Conditioned on the events $\AuxEvent(\gamcrit)$,
$\TailEvent(\gamcrit)$ and $\mathcal{C}(\gamcrit)$, and with our
choices of regularization parameter, we are guaranteed that the error
function $\DeltaHat$ satisfies all three of these constraints, and
hence that \mbox{$\DeltaHat \in \Gclass(\regpar, \regpartwo)$.}
Consequently, it suffices to establish a lower bound on
$\|g\|_\numobs$ that holds uniformly over the class $\Gclass(\regpar,
\regpartwo)$.  In particular, define the event
\begin{align}
\Bevent(\regpar, \regpartwo) & \defn \biggr \{
\|g\|_\numobs^2 \geq
\|g\|_2^2/2 \quad \mbox{for all $g \in \Gclass(\regpar, \regpartwo)$
  \mbox{where} $\|g\|_2 \geq \tildelcrit$} \biggr \}.
\end{align}
The following lemma shows that this event holds with high probability.
\blems
\label{LemLowerBound}
Under the conditions of Theorem~\ref{ThmPolyAchievable}, there are
universal constants $\plaincon_i$ such that
\begin{align}
\mprob[\Bevent(\regpar, \regpartwo)] & \geq 1 - \plaincon_1
\exp(-\plaincon_2 \numobs \gamcrit^2).
\end{align}
\elems
We note that this lemma can be interpreted as guaranteeting a version
of restricted strong convexity~\cite{Neg09} for the least-squares loss
function, suitably adapted to the non-parametric setting.  Since we do
not assume global boundedness, the proof of this lemma requires a
novel technical argument, one which combines a one-sided concentration
bound for non-negative random variables (Theorem 3.5 in Chung and
Lu~\cite{Chung06}) with the Sudakov minoration~\cite{Pisier89} for
Gaussian complexity.  We refer the reader to
Appendix~\ref{AppLemLowerBound} for all the details of the proof.\\

\noindent Using Lemma~\ref{LemLowerBound} and conditioning
$\Bevent(\regpar, \regpartwo)$, we are guaranteed that
$\|\DelHat\|_\numobs^2 \geq \|\DelHat\|_2^2/2$, and hence, combined
with our earlier bound~\eqref{EqnJunco}, we conclude that
\begin{align*}
\|\DeltaHat\|_2^2 \leq 4 \sqrt{\s}\regpar \|\DeltaHat\|_2 + 48
\sqrt{\s} \regpar \|\DeltaHat\|_2 + 24 \s \regpar^2 + 50 \s
\regpartwo.
\end{align*}
Hence $\|\DeltaHat\|_\numobs \leq \|\DeltaHat\|_2 \leq C
\max(\sqrt{\s}\regpar, \sqrt{\s \regpartwo})$, completing the proof of
the claim in the third case. \\

In summary, the entire proof is based on conditioning on the three
events $\TailEvent(\gamcrit)$, $\AuxEvent(\regpar)$ and
$\Bevent(\regpar, \regpartwo)$.  From the bound~\eqref{EqnAuxBound} as
well as Lemmas~\ref{LemNetherlands} and~\ref{LemLowerBound}, we have
\begin{align*}
\mprob \big[ \TailEvent(\gamcrit) \cap \AuxEvent(\regpar) \cap
  \Bevent(\regpar, \regpartwo)\cap \mathcal{C}(\gamcrit) \big] & \geq
1- c_1 \exp \big( - c_2 \numobs \gamcrit^2 \big),
\end{align*}
thereby showing that $\max \{ \|\fhat - \fstar\|_\numobs^2, \|\fhat -
\fstar\|_2^2 \} \leq C \max(\s \regpar^2, \s \regpartwo)$ with the
claimed probability. This completes the proof of
Theorem~\ref{ThmPolyAchievable}.

\subsection{Proof of Theorem~\ref{ThmLower}}
\label{SecLowerBound}

We now turn to the proof of the minimax lower bounds stated in
Theorem~\ref{ThmLower}.  For both parts (a) and (b), the first step is
to follow a standard reduction to testing
(e.g.,~\cite{Hasminskii78,YanBar99,Yu}) so as to obtain a lower bound
on the minimax error $\MiniMaxPop$ in terms of the probability of
error in a multi-way hypothesis testing.  We then apply different
forms of the Fano inequality~\cite{YanBar99,Yu} in order to lower
bound the probability of error in this testing problem.  Obtaining
useful bounds requires a precise characterization of the metric
entropy structure of $\MyBigClass$, as stated in
Lemma~\ref{LemMetEnt}.

\subsubsection{Reduction to testing}

We begin with the reduction to a testing problem.  Let $\{f^1, \ldots,
f^M \}$ be a $\epspack$-packing of $\Fclass$ in the
$\|\cdot\|_2$-norm, and let $\PackVar$ be a random variable uniformly
distributed over the index set \mbox{$[M] \defn \{1, 2, \ldots,
M\}$.}  Note that we are using $M$ as a shorthand for the
packing number $M(\epspack; \Fclass, \|\cdot\|_2)$.  A standard
argument (e.g.,~\cite{Hasminskii78,YanBar99,Yu}) then yields the lower
bound
\begin{align}
\label{ExpectBound}
\inf_{\fhat} \sup_{\fstar \in \Fclass} \mprob \big[ \|\fhat -
\fstar\|_2^2 \geq \epspack^2/2 \big] & \geq \inf_{\PackVarHat}
\mprob[\PackVarHat \neq \PackVar],
\end{align}
where the infimum on the right-hand side is taken over all estimators
$\PackVarHat$ that are measurable functions of the data, and take
values in the index set $[M]$.

Note that $\mprob[\PackVarHat \neq \PackVar]$ corresponds to the error
probability in a multi-way hypothesis test, where the probability is
taken over the random choice of $\PackVar$, the randomness of the
design points $\Xdes \defn \{\exi\}_{i=1}^\numobs$, and the randomness
of the observations $\Ysam \defn \{\yi\}_{i=1}^\numobs$.  Our initial
analysis is performed conditionally on the design points, so that the
only remaining randomness in the observations $\Ysam$ comes from the
observation noise $\{\wi\}_{i=1}^\numobs$.  From Fano's
inequality~\cite{Cover}, for any estimator $\PackVarHat$, we have
$\mprob \big[\PackVarHat \neq \PackVar \, \mid \, \Xdes \big ] \geq 1
- \frac{\InfoX(\PackVar; \Ysam) + \log 2}{\log M}$, where
$\InfoX(\PackVar; \Ysam)$ denotes the mutual information between
$\PackVar$ and $\Ysam$ with $\Xdes$ fixed.  Taking expectations over
$\Xdes$, we obtain the lower bound
\begin{align}
\label{EqnFanoAve}
\mprob \big[\PackVarHat \neq \PackVar \big ] & \geq 1 -
\frac{\Exs_{\Xdes}\big[\InfoX(\PackVar; \Ysam) \big] + \log 2}{\log
M}.
\end{align}
The remainder of the proof consists of constructing appropriate
packing sets of $\Fclass$, and obtaining good upper bounds on the
mutual information term in the lower bound~\eqref{EqnFanoAve}.

\subsubsection{Constructing appropriate packings}

We begin with results on packing numbers. Recall that $\log M(\delta;
\Fclass, \|\cdot\|_2)$ denotes the $\delta$-packing entropy of
$\Fclass$ in the $\|\cdot\|_2$ norm.
\blems
\label{LemMetEnt}

\begin{enumerate}
\item[(a)] For all $\delta \in (0,1)$ and $\spindex \leq \pdim/4$, we
have
\begin{align}
\label{EqnCovering}
\log M(\delta; \Fclass, \|\cdot\|_2) & = \mathcal{O} \big( \spindex \; \log
M(\frac{\delta}{\sqrt{\kdim}}; \Ball_\Hil(1), \|\cdot\|_2) + \spindex
\log \frac{\pdim}{\spindex} \big).
\end{align}
\item[(b)] For a Hilbert class with logarithmic metric
entropy~\eqref{EqnLogMet} and such that $\|f\|_2 \leq \|f\|_\Hil$,
there exists set $\{f^1, \ldots, f^M\}$ with
$\log M \; \geq \; \UniCon \, \big \{ \spindex
\log(\pdim/\spindex) + \spindex \kerrank \big \}$, and
\begin{equation}
\label{EqnLogMetPacking}
\delta \; \leq \; \|f^k - f^m\|_2 \; \leq \; 8 \delta \qquad \mbox{for
all $k \neq m \in \{1,2, \ldots, M\}$.}
\end{equation}
\end{enumerate}
\elems
\noindent The proof, provided in Appendix~\ref{AppLemMetEnt}, is
combinatorial in nature.  We now turn to the proofs of parts (a) and
(b) of Theorem~\ref{ThmLower}.


\subsubsection{Proof of Theorem~\ref{ThmLower}(a)}

In order to prove this claim, it remains to exploit
Lemma~\ref{LemMetEnt} in an appropriate way, and to upper bound the
resulting mutual information.  For the latter step, we make use of the
generalized Fano approach (e.g.,~\cite{Yu}).

From Lemma~\ref{LemMetEnt}, we can find a set $\{f^1, \ldots,
f^\GoodPackSize\}$ that is a $\delta$-packing of $\Fclass$ in
$\ell_2$-norm, and such that $\|f^k - f^\ell\|_2 \leq 8 \delta$ for
all $k, \ell \in [\GoodPackSize]$.  For $k = 1, \ldots,
\GoodPackSize$, let $\Qprob^k$ denote the conditional distribution of
$\Ysam$ conditioned on $\Xdes$ and the event $\{\PackVar = k\}$, and
let $\kull{\Qprob^k}{\Qprob^\ell}$ denote the Kullback-Leibler
divergence.  From the convexity of mutual information~\cite{Cover}, we
have the upper bound $\InfoX(\PackVar; \Ysam) \leq
\frac{1}{{\GoodPackSize \choose 2}} \sum_{k, \ell = 1}^\GoodPackSize
\kull{\Qprob^k}{\Qprob^\ell}$.  Given our linear observation
model~\eqref{EqnLinObs}, we have
\begin{align*}
\kull{\Qprob^k}{\Qprob^\ell} & = \frac{1}{2 \sigma^2} \sum_{i=1}^\numobs
\big(f^k(\xsam{i}) - f^\ell(\xsam{i}) \big)^2 \; = \; \frac{\numobs \;
\|\f^k - f^\ell\|_\numobs^2}{2 },
\end{align*}
and hence
\begin{align*}
\Exs_\Xdes \big[ \InfoX(\Ysam; \PackVar) \big] & \leq
  \frac{\numobs}{2} \frac{1}{{\GoodPackSize \choose 2}} \sum_{k, \ell =
  1}^\GoodPackSize \Exs_\Xdes[\|f^k - f^\ell\|_\numobs^2] \; = \;
  \frac{\numobs}{2} \frac{1}{{\GoodPackSize \choose 2}} \sum_{k, \ell =
  1}^\GoodPackSize \|f^k - f^\ell\|_2^2.
\end{align*}
Since our packing satisfies $\|f^k - f^\ell\|_2^2 \leq 64 \delta^2$, we
conclude that
\begin{align*}
\Exs_\Xdes \big[ \InfoX(\Ysam; \PackVar) \big] & \leq 32 \numobs
\delta^2.
\end{align*}
From the Fano bound~\eqref{EqnFanoAve}, for any $\delta > 0$ such that
$\frac{32 \numobs \delta^2 + \log 2}{\log \GoodPackSize} <
\frac{1}{4}$, then we are guaranteed that $\mprob[\PackVarHat \neq
\PackVar] \geq \frac{3}{4}$.  From Lemma~\ref{LemMetEnt}(b), our
packing set satisfies $\log \GoodPackSize \geq \UniCon \big \{
\spindex \mdim + \spindex \log (\pdim/\spindex) \big \}$, so that so
that the choice $\delta^2 = \UniCon' \, \big \{ \frac{\spindex
\mdim}{\numobs} + \frac{\spindex \log(\pdim/\spindex)}{\numobs} \big
\}$, for a suitably small $\UniCon' > 0$, can be used to guarantee the
error bound $\mprob[\PackVarHat \neq \PackVar] \geq \frac{3}{4}$.

\subsubsection{Proof of Theorem~\ref{ThmLower}(b)}
\label{SecProofThmLower}
In this case, we use an upper bounding technique due to Yang and
Barron~\cite{YanBar99} in order to upper bound the mutual
information. Although the argument is essentially the same, it does
not follow verbatim from their claims---in particular, there are some
slight differences due to our initial conditioning---so that we
provide the details here.  By definition of the mutual information, we
have
\begin{align*}
I_\Xdes(\PackVar; \Ysam) & = \frac{1}{\GoodPackSize}
\sum_{k=1}^\GoodPackSize \kull{\Qprob^k}{\mprob_Y},
\end{align*}
where $\Qprob^k$ denotes the conditional distribution of $\Ysam$ given
$\PackVar = k$ and still with $\Xdes$ fixed, whereas $\mprob_Y$
denotes the marginal distribution of $\mprob_Y$. 

First we define \emph{covering numbers}. Let $(\GenSet, \rho)$ be a
totally bounded metric space, consisting of a set $\GenSet$ and a
metric $\rho: \GenSet \times \GenSet \rightarrow \mathbb{R}_+$.  An
$\epsilon$-covering set of $\GenSet$ is a collection $\{\f^1, \ldots,
\f^N \}$ of functions such that for all $f \in \GenSet$ there exists
$k \in \{ 1,2,...,N\}$ such that $\rho(f, f^k) \leq \epsilon$. The
$\epsilon$-covering number $N(\epsilon; \GenSet, \rho)$ is the
cardinality of the smallest $\epsilon$-covering set.

Now let $\{g^1,\ldots, g^N\}$ be an $\epsilon$-cover of $\Fclass$ in the
$\|\cdot\|_2$ norm, for a tolerance $\epsilon$ to be chosen.  As
argued in Yang and Barron~\cite{YanBar99}, we have
\begin{align*}
I_\Xdes(\PackVar; \Ysam) \; = \; \frac{1}{\GoodPackSize}
\sum_{j=1}^\GoodPackSize \kull{\Qprob^j}{\mprob_Y} & \leq
\kull{\Qprob^k}{\frac{1}{N} \sum_{k=1}^N \mprob^k},
\end{align*}
where $\mprob^\ell$ denotes the conditional distribution of $\Ysam$ given
$g^\ell$ and $\Xdes$.  For each $\ell$, let us choose \mbox{$\ell^*(k) \in \arg
\min_{\ell = 1, \ldots, N} \|g^\ell - f^k\|_2$.}  We then have the upper
bound
\begin{align*}
I_\Xdes(\PackVar; \Ysam) & \leq \frac{1}{M} \sum_{k=1}^M \big \{ \log
N + \frac{\numobs}{2} \|g^{\ell^*(k)} - f^k\|_\numobs^2 \big \}.
\end{align*}
Taking expectations over $\Xdes$, we obtain
\begin{align*}
\Exs_\Xdes[\InfoX(\PackVar; \Ysam)] & \leq \frac{1}{M} \sum_{k=1}^M
\big \{ \log N + \frac{\numobs}{2} \Exs_\Xdes [\|g^{\ell^*(k)} -
f^k\|_\numobs^2] \big \} \\
& \leq \log N + \frac{\numobs}{2} \; \epsilon^2,
\end{align*}
where the final inequality follows from the choice of our covering
set.

From this point, we can follow the same steps as Yang and
Barron~\cite{YanBar99}.  The polynomial scaling~\eqref{EqnPolyMet} of
the metric entropy guarantees that their conditions are satisfied, and
we conclude that the minimax error is lower bounded any $\epspack > 0$
such that $\numobs \epspack^2 \; \geq \; C \log N(\epspack; \Fclass,
\|\cdot\|_2)$.  From Lemma~\ref{LemMetEnt} and the assumed
scaling~\eqref{EqnPolyMet}, it is equivalent to solve the equation
\begin{align*}
\numobs \epspack^2 & \; \geq \; C \biggr \{ \spindex \log(\pdim/\spindex) +
\spindex (\sqrt{\spindex}/\epspack)^{1/\smooth} \biggr \},
\end{align*}
from which some algebra yields $\epspack^2 = \UniCon \big \{
\frac{\spindex \log(\pdim/\spindex)}{\numobs} + \spindex \big(
\frac{1}{\numobs} \big)^{\frac{2 \smooth}{2 \smooth + 1}} \big\}$ as a
suitable choice.


\subsection{Proof of Theorem~\ref{ThmFast}}
\label{SecProofFast}

Recall the definition of $\MyNewClass$ and $\MyNewSmallClass$ from
Section~\ref{SecKYCompare} which ensures $f^*$ is uniformly bounded by
$B$. In order to establish upper bounds on the minimax rate in
$\LTP$-error over $\MyNewClass$, we we analyze a least-squares
estimator (not the same as the original M-estimator~\eqref{EqnMest})
constrained to $\MyNewClass$:
\begin{align}
\label{EqnEstimate}
\fhat & \in \arg \min_{f \in \MyNewClass} \sum_{i=1}^{n}{(y_i - \ybar
  - f(\exi))^2}.
\end{align}
Since our goal is to upper bound the minimax rate in $\LTP$ error, it
is sufficient to upper bound the $\LTP$-norm of $\fhat-\fstar$ where
$\fhat$ is any solution to ~\eqref{EqnEstimate}. The proof shares many
steps with the proof of Theorem~\ref{ThmPolyAchievable}.  First, the
same reasoning shows that the error $\DelHat \defn \fhat - \fstar$
satisfies the basic inequality
\begin{align}
\label{EqGeneral}	
\frac{1}{\numobs} \sum_{i=1}^{n}{\DelHat^2(x_i)} & \leq
\frac{2}{\numobs} |\sum_{i=1}^{n}{w_i\DelHat(x_i)}|+ |\ybar-\mu | \big
| \frac{1}{\numobs} \sum_{i=1}^\numobs {\DeltaHat(\exi)} \big |.
\end{align}
Recall the definition~\eqref{EqnDefnSmallRate} of the critical rate
$\UnifBoundRate$. Once again, we first control the term error due to
estimating the mean $|\ybar-\mu|$. Noting that $\ybar-\mu =
\frac{1}{n}\sum_{i=1}^{n}(y_i - \mu)$ and the random variable $y_i -
\mu = \fstar(\exi) + \wi$ is sub-Gaussian with parameter
$\sqrt{B+1}$. This follows since $f^*$ is bounded by $B$ and using
standard results on sums of independent sub-Gaussian random variables
(see e.g. Lemma 1.4 in Buldygin and
Kozachenko~\cite{BulKoz}). Therefore
\begin{equation*}
\mP(|\ybar-\mu| > t) \leq 2 \exp(- \frac{n t^2}{2 (B+1)}).	
\end{equation*}	
Setting $\mathcal{A}(\UnifBoundRate) = \{ |\ybar-\mu| \leq
\sqrt{B}\UnifBoundRate \}$, it is clear that
\begin{equation*}
\mP[\mathcal{A}(\UnifBoundRate)] \geq 1 - 2 \exp(- \frac{n
  \UnifBoundRate^2}{4}).
\end{equation*}	
For the remainder of the proof, we condition on the event
$\mathcal{A}(\UnifBoundRate)$, in which case
equation~\eqref{EqRiskBound} simplifies to
\begin{equation}
\label{EqnBoundBasicTwo}
\frac{1}{2}\|\DeltaHat\|_\numobs^2 \leq \big | \frac{1}{\numobs}
\sum_{i=1}^\numobs {\wi\DeltaHat(\exi)} \big | +
\sqrt{B}\UnifBoundRate\|\DeltaHat\|_n.
\end{equation}
Here we have used the fact that $\big | \frac{1}{\numobs}
\sum_{i=1}^\numobs {\DeltaHat(\exi)} \big | \leq
\|\DeltaHat\|_\numobs$, by the Cauchy-Schwartz inequality.

Now we control the Gaussian complexity term $\big | \frac{1}{\numobs}
\sum_{i=1}^\numobs {\wi\DeltaHat(\exi)} \big |$. For any fixed subset
$S$, define the random variable
\begin{align}
\label{GaussCompBound}
\RvarHat(w,t; \HACK) & \defn \sup_ { \substack{\Delta \in \HACK
    \\ \|\Delta\|_\numobs \leq t}} \big| \frac{1}{\numobs}
\sum_{i=1}^\numobs w_i \Delta(x_i) \big|.
\end{align}
We first bound this random variable for a fixed subset $S$ of size $2
\s$, and then take the union bound over all ${\pdim \choose 2 \s}$
possible subsets.
\blems
\label{LemSquash}
Assume that the RKHS $\Hil$ has eigenvalues $(\mu_k)_{k=1}^{\infty}$
that satisfy $\mu_k \simeq k^{-2 \alpha}$ and eigenfunctions such that
$\|\phi_k\|_{\infty} \leq C$.  Then we have
\begin{align}
\mprob \big[ \exists t > 0 \mbox{ such that } \RvarHat(w,t; \HACK)
  \geq 16 B \, C \,\sqrt{\frac{\s^{1/\alpha}\log \s}{\numobs}} + 3 t
  \UnifBoundRate \big] & \leq c_1 \exp(-9 \numobs \UnifBoundRate^2).
\end{align}
\elems
\noindent The proof of Lemma~\ref{LemSquash} is provided
Appendix~\ref{AppLemRadGauss}. Returning to
inequality~\eqref{EqnBoundBasicTwo}, we note that by definition,
\begin{align*}
\frac{2}{\numobs} |\sum_{i=1}^{n}{w_i\DelHat(x_i)}| & \leq \max_{|S| =
  2 \s} \RvarHat(w,\|\DelHat\|_\numobs; \HACK).
\end{align*}
Lemma~\ref{LemSquash} combined with the union bound implies that
\begin{align*}
\max_{|S| = 2 \s} \RvarHat(w,\|\DelHat\|_\numobs; \HACK) & \leq 16 B
\, C \,\sqrt{\frac{\s^{1/\alpha}\log \s}{\numobs}} + 3 \UnifBoundRate
\|\DelHat\|_\numobs
\end{align*}
with probability at least $1 - c_1 {\pdim \choose 2 \s} \; \exp(-3
\numobs \UnifBoundRate^2)$.  Our choice~\eqref{EqnDefnSmallRate} of
$\UnifBoundRate$ ensures that this probability is at least $1 - c_1
\exp(-c_2 \numobs \UnifBoundRate^2)$.  Combined with the basic
inequality~\eqref{EqnBoundBasicTwo}, we conclude that
\begin{align}
\label{EqnHilary}
\|\DeltaHat\|_n^2 & \leq 32 B \, C \,\sqrt{\frac{\s^{1/\alpha}\log
    \s}{\numobs}} + 7 B \UnifBoundRate \, \|\DeltaHat\|_\numobs
\end{align}
with probability $1 - c_1 \exp(-c_2 \numobs \UnifBoundRate^2)$.  

By definition~\eqref{EqnDefnSmallRate} of $\UnifBoundRate$, the
bound~\eqref{EqnHilary} implies that $\|\DeltaHat\|_\numobs =
\order(\UnifBoundRate)$ with high probability.  In order to translate
this claim into a bound on $\|\DelHat\|_2$, we require the following
result:
\blems
\label{LemGeneralBound}
There exist universal constants $(\plaincon, \plaincon_1, \plaincon_2)$
such that for all $t \geq \plaincon \UnifBoundRate$, we have
\begin{equation}
\label{EqnSandwichBound}
\frac{\|g\|_2}{2} \; \leq \; \|g\|_\numobs \; \leq \; \frac{3}{2}
\|g\|_2 \qquad \mbox{for all $g \in \HACK$ with $\|g\|_2 \geq t$}
\end{equation}
with probability at least $1 -\plaincon_1 \exp(-\plaincon_2 \numobs
t^2)$.
\elems 

\spro
The bound~\eqref{EqnSandwichBound} follows by applying
Lemma~\ref{LemGeneral} in Appendix~\ref{AppLemEquivNorms} with
$\Gclass = \HACK$ and $\smallbou = 2 \GLOB$.  The critical radius from
equation~\eqref{EqnDefnEpscrit} needs to satisfy the relation
$\LocGauss(\epscrit; \HACK) \leq \frac{\epscrit^2}{40}$.  From
Lemma~\ref{LemSmallGauss}, the choice $\epscrit^2 = 320 \GLOB \, C
\,\sqrt{\frac{\s^{1/\alpha}\log \s}{n}}$ satisfies this relation.  By
definition~\eqref{EqnDefnSmallRate} of $\UnifBoundRate$, we have
$\UnifBoundRate \geq c \epscrit$ for some universal constant $c$,
which completes the proof.  \fpro

This lemma implies that with probability at least $1 -\plaincon_1
\exp(-\plaincon_2 B \numobs \UnifBoundRate^2)$, we have $\|\DelHat\|_2
\leq 2 \|\DelHat\|_n + C \UnifBoundRate$.  Combined with our earlier
upper bound on $\|\DelHat\|_\numobs$, this completes the proof of
Theorem~\ref{ThmFast}.


\section{Discussion}
\label{SecDiscuss}

In this paper, we have studied estimation in the class of sparse
additive models defined by univariate reproducing kernel Hilbert
spaces.  In conjunction, Theorems~\ref{ThmPolyAchievable}
and~\ref{ThmLower} provide a precise characterization of the
minimax-optimal rates for estimating $\fstar$ in the $\LTP$-norm for
various kernel classes with bounded univariate functions. These
classes include finite-rank kernels (with logarithmic metric entropy),
as well as kernels with polynomially decaying eigenvalues (and hence
polynomial metric entropy).  In order to establish achievable rates,
we analyzed a simple $M$-estimator based on regularizing the
least-squares loss with two kinds of $\ell_1$-based norms, one defined
by the univariate Hilbert norm and the other by the univariate
empirical norm.  On the other hand, we obtained our lower bounds by a
combination of approximation-theoretic and information-theoretic
techniques.

An important feature of our analysis is we assume only that each
univariate function is bounded, but do not assume that the
multivariate function class is bounded.  As discussed in
Section~\ref{SecKYCompare}, imposing a global boundedness condition in
the high-dimensional setting could lead to a substantially smaller
function classes; for instance, for Sobolev classes and sparsity $\s =
\Omega(\sqrt{\numobs})$, Theorem~\ref{ThmFast} shows that it is
possible to obtain much faster rates than the optimal rates for the
class of sparse additive models with univariate functions
bounded. Theorem~\ref{ThmFast} in our paper shows that the rates
obtained under global boundedness conditions are not minimax optimal
for Sobolev spaces in the regime $\s = \Omega(\sqrt{\numobs})$.

There are a number of ways in which this work could be extended. For
instance, although our analysis was based on assuming independence of
the covariates $x_j$, $j=1,2, \ldots \pdim$, it would be interesting
to investigate the case when the random variables are endowed with
some correlation structure. One might expect some changes in the
optimal rates, particularly if many of the variables are strongly
dependent.  This work considered only the function class consisting of
sums of univariate functions; a natural extension would be to consider
nested non-parametric classes formed of sums over hierarchies of
subsets of variables.  Analysis in this case would require dealing
with dependencies between the different functions and is left for
future research.

\subsection*{Acknowledgements}
This work was partially supported by NSF grants DMS-0605165 and
DMS-0907632 to MJW and BY.  In addition, BY was partially supported by
the NSF grant SES-0835531 (CDI), the SRO grant (INSERT NUMBER) and the
Purdue grant (INSERT NUMBER). MJW was also partially supported AFOSR
Grant FA9550-09-1-0466.  During this work, GR was financially
supported by a Berkeley Graduate Fellowship.


\appendix

\section{A general result on equivalence of $\LTP$ and $\LTPn$ norms}
\label{AppLemEquivNorms}

Since it is required in a number of our proofs, we begin by stating
and proving a general result that provides uniform control on the
difference between the empirical $\|\cdot\|_\numobs$ and population
$\|\cdot\|_2$ norms over a uniformly bounded function class
$\FUNCLASS$.  We impose two conditions on this class:
\begin{enumerate}
\item[(a)] it is uniformly bounded, meaning that there is some
  $\smallbou \geq 1$ such that $\|g\|_\infty \leq \smallbou$ for all
  $g \in \FUNCLASS$.
\item[(b)] it is star-shaped, meaning that if $g \in \FUNCLASS$, then
  $\lambda g \in \FUNCLASS$ for all $\lambda \in [0,1]$.
\end{enumerate}
For each co-ordinate, the Hilbert ball $\Ball_\Hil(2)$ satisfies both
of these conditions; we use $\FUNCLASS = \Ball_\Hil(2)$. (To be clear,
we cannot apply this result to the multivariate function class
$\MyBigClass$, since it is not uniformly bounded.)

Let $\{\rade{i}\}_{i=1}^\numobs$ be an i.i.d. sequence of Rademacher
variables, and let $\{x_i\}_{i=1}^\numobs$ be an i.i.d. sequence of
variables from $\Xcal$, drawn according to some distribution $\Qprob$.
For each $t > 0$, we define the local Rademacher complexity
\begin{equation}
\label{EqnDefnQrad}
\LocRad(t, \FUNCLASS) \; \defn \Exs_{x, \radepl}
\big[\sup_{\substack{\|g\|_2 \leq t \\ g \in \FUNCLASS}}
  \frac{1}{\numobs} \sum_{i=1}^\numobs \rade{i} g(\exi) \big]
\end{equation}
We let $\epscrit$ denote the smallest solution (of size at least
$1/\sqrt{\numobs}$) to the inequality
\begin{equation}
\label{EqnDefnEpscrit}
\LocRad(\epscrit, \FUNCLASS) = \frac{\epscrit^2}{40},
\end{equation}
where our scaling by the constant $40$ is for later theoretical
convenience.  Such an $\epscrit$ exists, because the star-shaped
property implies that the function $\LocRad(t, \FUNCLASS)/t$ is
non-increasing in $t$.  This quantity corresponds to the critical rate
associated with the population Rademacher complexity.  For any $t \geq
\epscrit$, we define the event $\Event(t) \defn \big \{
\sup_{\substack{g \in \FUNCLASS \\ \|g\|_2 \leq t}} \big|
\|g\|_\numobs - \|g\|_2 \big| \geq \frac{\smallbou t}{2} \big \}$.
\blems
\label{LemGeneral}
 Suppose that $\|g\|_\infty \leq \smallbou$ for all $g \in
 \FUNCLASS$. Then there exist universal constants $(\plaincon_1,
 \plaincon_2)$ such that for any $t \geq \epscrit$,
\begin{align}
\label{EqnGeneral}
\mprob \big[ \Event(t) \big] & \leq \plaincon_1 \exp(-\plaincon_2
  \numobs t^2).
\end{align}
In addition, for any $g \in \FUNCLASS$ with $\|g\|_2 \geq t$, we have
$\|g\|_\numobs \; \leq \; \|g\|_2(1 + \frac{b}{2})$, and moreover, for
all $g \in \FUNCLASS$ with $\|g\|_2 \geq bt$, we have
\begin{equation}
\label{EqnSandwich}
\frac{1}{2}\|g\|_2\; \leq \;\|g\|_\numobs \; \leq \; \frac{3}{2}\|g\|_2,
\end{equation}
both with probability at least $1 -\plaincon_1 \exp(-\plaincon_2
\numobs t^2)$.
\elems
\noindent Lemma~\ref{LemGeneral} follows from a relatively
straightforward adaptation of known results (e.g., Lemma 5.16 in van
de Geer~\cite{vandeGeer} and Theorem 2.1 in Bartlett et
al.~\cite{Bartlett05}), so we omit the proof details here.

\section{Proof of Lemma~\ref{LemNetherlands}}
\label{AppLemNetherlands}

The proof of this lemma is based on peeling and weighting techniques
from empirical theory~\cite{Alex87, vandeGeer} combined with results
on the local Rademacher and Gaussian complexities~\cite{Bartlett05,
  Mendelson02}.  For each univariate Hilbert space $\Hil_j=\Hil$, let
us introduce the random variables
\begin{equation}
\label{EqnDefnRvar}
\RvarHat(w, t; \Hil) \defn \sup_{\substack{\|g_j\|_{\Hil} \leq 1
    \\ \|g_j\|_\numobs \leq t}} \big| \frac{1}{\numobs}
\sum_{i=1}^\numobs w_i g_j(x_{ij}) \big|, \quad \mbox{ and } \quad \Rvar(w,
t; \Hil) \defn \Exs_x \biggr[ \sup_{\substack{\|g_j\|_\Hil \leq 1
      \\ \|g_j\|_2 \leq t}} \big| \frac{1}{\numobs} \sum_{i=1}^\numobs
  w_i g_j(x_{ij}) \big| \biggr],
\end{equation}
where $w_i \sim N(0,1)$ are i.i.d. standard normal.  The empirical and
population Gaussian complexities are given by
\begin{align}
\label{EqnDefnRvarExp}
\LocGaussEmp(t, \Hil) \; \defn \Exs_{w} \big[\RvarHat(w; t, \Hil)
  \big] \quad \mbox{ and } \quad
\LocGauss(t, \Hil) \; \defn \Exs_{w} \big[\Rvar(w; t, \Hil) \big].
\end{align}
For future reference, we note that in the case of a univariate Hilbert
space $\Hil$ with eigenvalues $\{\mu_k\}_{k=1}^\infty$, results in
Mendelson~\cite{Mendelson02} imply that there are universal constants
$c_\ell \leq c_u$ such that for all $t^2 \geq 1/\numobs$, we have
\begin{equation}
\label{EqnMend}
\frac{c_\ell}{\sqrt{\numobs}} \big[ \sum_{k=1}^\infty \min \{ t^2,
  \mu_k \}\big]^{1/2} \; \leq \; \LocGauss(t, \Hil) \; \leq \;
\frac{c_u}{\sqrt{\numobs}} \big[ \sum_{k=1}^\infty \min \{ t^2, \mu_k
  \}\big]^{1/2},
\end{equation}
for all $j$
The same bounds also hold for the local Rademacher complexities for
Reproducing kernel Hilbert spaces.

Let $\delunihat{j} > 0$ denote the smallest positive solution $r$ of the
inequality 
\begin{align}
\label{EqnDefnDelHatUni}
\LocGaussEmp(r, \Hil) & \leq 4 \, r^2.
\end{align}
The function $\LocGaussEmp(r, \Hil)$ defines the local Gaussian
complexity of the kernel class in co-ordinate $j$.  Recall the
bounds~\eqref{EqnMend} that apply to both the empirical and population
Gaussian complexities.  Recall that the critical univariate rate
$\unicrit$ is defined in terms of the population Gaussian complexity
(see equation~\eqref{EqnCritCon}).

\subsection{Some auxiliary results} 

In order to prove Lemma~\ref{LemNetherlands}, we also need some
auxiliary results, stated below as Lemmas~\ref{LemGconc}
and~\ref{LemKerEquiv}.

\blems
\label{LemGconc}
For any function class $\FUNCLASS$ and all $\delta \geq 0$, we have
\begin{subequations}
\begin{align}
\label{EqnLedOne}
\mprob \big [|\RvarHat(w,t,\FUNCLASS) - \LocGaussEmp(t,\FUNCLASS)| \geq \delta t
  \big ] & \leq 2 \exp \big(-\frac{n \delta^2}{2} \big), \quad
\mbox{and} \\
\label{EqnLedTwo}	
\mprob \big[ |\Rvar(w,t,\FUNCLASS) - \LocGauss(t,\FUNCLASS)| \geq \delta t \big
] & \leq 2 \exp\biggr(-\frac{\numobs \delta^2}{2} \biggr).
\end{align}
\end{subequations}
\elems

\spro
We have
\begin{align*}
|\RvarHat(w,t,\FUNCLASS) - \RvarHat(w',t,\FUNCLASS)| \; \leq \;
\sup_{\substack{g \in \FUNCLASS \\ \|g\|_\numobs \leq t}}
\frac{1}{n}|\sum_{i=1}^{\numobs}(w_i - w'_i)g(x_i)| \; \leq \;
\frac{t}{\sqrt{n}} \|w-w'\|_2,
\end{align*}
showing that $\RvarHat(w,t,\FUNCLASS)$ is $\frac{t}{\sqrt{n}}$-Lipschitz
with respect to the $\ell_2$ norm.  Consequently, concentration for
Lipschitz functions of Gaussian random variables~\cite{Ledoux01}
yields the tail bound~\eqref{EqnLedOne}.  Turning to the quantity
$\Rvar(w,t,\Hil)$, a similar argument yields that
\begin{align*}
|\Rvar(w,t,\FUNCLASS) - \Rvar(w',t,\FUNCLASS)| & \leq \Exs_x
\big[\sup_{\substack{g \in \FUNCLASS \\ \|g\|_2 \leq t}}
  \frac{1}{n}|\sum_{i=1}^{\numobs}(w_i - w'_i)g(x_i)| \big] \\
& \leq \; \; \sup_{\substack{g \in \FUNCLASS \\ \|g\|_2 \leq t}}
\Exs_x \big[ \big(\frac{1}{\numobs} \sum_{i=1}^\numobs g^2(x_i))^{1/2}
  \big] \; \|w-w'\|_2 \; \leq \; \frac{t}{\sqrt{\numobs}} \, \|w -
w'\|_2,
\end{align*}
where the final step uses Jensen's inequality and the fact that
$\Exs_x[g^2(x_i)] \leq t^2$ for all $i = 1, \ldots, \numobs$.
The same reasoning then yields the tail bound~\eqref{EqnLedTwo}.
\fpro

Our second lemma involves the event $\AnnoyEvent(\gamcrit) \defn \big
\{ \delhatuni{j} \leq \gamcrit, \quad \mbox{for all $j = 1, 2, \ldots,
  \pdim$} \big \}$, where we recall the
definition~\eqref{EqnDefnDelHatUni} of $\delhatuni{j}$, and that
$\gamcrit \defn \kappa \max \big \{ \unicrit, \sqrt{\frac{\log
    \pdim}{\numobs}} \big \}$.
\blems
\label{LemKerEquiv}
For all $1 \leq j \leq \pdim$, we have
\begin{align}
\mprob \big[ \delhatuni{j} \leq \gamcrit \big] & \geq 1 - \plaincon_1
\exp(-\plaincon_2 n \gamcrit^2).
\end{align}
\elems 

\spro

We first bound the probability of the event $\{\delhatuni{j} >
\gamcrit \}$ for a fixed $\Hil_j$.  Let $g \in \Ball_{\Hil_j}(1)$ be
any function such that $\|g\|_2 > t \geq \unicrit$.  Then conditioned
on the sandwich relation~\eqref{EqnSandwich} with $\smallbou = 1$, we
are guaranteed that $\|g\|_\numobs > \frac{t}{2}$.  Taking the
contrapositive, we conclude that $\|g\|_\numobs \leq \frac{t}{2}$
implies $\|g\|_2 \leq t$, and hence that $\RvarHat(w,t/2,\Hil) \leq
\Rvar(w,t,\Hil)$ for all $t \geq \unicrit$, under the stated
conditioning.

For any $t \geq \unicrit$, the inequalities~\eqref{EqnSandwich},
~\eqref{EqnLedOne} and ~\eqref{EqnLedTwo} hold with probability at
least \mbox{$1-\plaincon_1 \exp(-\plaincon_2 \numobs t^2)$.}
Conditioning on these inequalities, we can set $t = \gamcrit >
\unicrit$, and thereby obtain
\begin{align*}
\LocGaussEmp(\gamcrit,\Hil) & \stackrel{(a)}{\leq}
\RvarHat(w,\gamcrit,\Hil) + \gamcrit^2 \\
& \stackrel{(b)}{\leq} \; \Rvar(w,2 \gamcrit,\Hil) + \gamcrit^2 \\
& \stackrel{(c)}{\leq} \; \LocGauss(2 \gamcrit,\Hil) + 2 \gamcrit^2 \\
& \stackrel{(d)}{\leq} 4 \gamcrit^2,
\end{align*}
where inequality (a) follows from the bound~\eqref{EqnLedOne},
inequality (b) follows the initial argument, inequality (c) follows
from the bound~\eqref{EqnLedTwo}, and inequality (d) follows since $2
\gamcrit > \epscrit$ and the definition of $\epscrit$.  

By the definition of $\delhatuni{j}$ as the minimal $t$ such that
$\LocGaussEmp(t,\Hil) \leq 4 t^2$, we conclude that for each fixed
$j = 1, \ldots, \numobs$, we have $\delhatuni{j} \leq \gamcrit$ with
probability at least \mbox{$1-\plaincon_1 \exp(-\plaincon_2 \numobs
  \gamcrit^2)$.}  Finally, the uniformity over $j = 1, 2, \ldots,
\pdim$ follows from the union bound and our choice of $\gamcrit \geq
\kappa \sqrt{\frac{\log \pdim}{\numobs}}$.

\fpro

\subsection{Main argument to prove Lemma~\ref{LemNetherlands}}

We can now proceed with the proof of Lemma~\ref{LemNetherlands}.
Combining Lemma~\ref{LemKerEquiv} with the union bound over $j=1,2,
\ldots,\pdim$, we conclude that that
\begin{align*}
\mprob[\AnnoyEvent(\gamcrit)] \geq 1 - \plaincon_1 \exp(-\plaincon_2
\numobs \gamcrit^2),
\end{align*}
as long as $c_2 \geq 1$. For the remainder of our proofs, we condition
on the event $\AnnoyEvent(\gamcrit)$.  In particular, our goal is to
prove that
\begin{align}
\label{EqnNetherlandsRandom}
\big| \frac{1}{\numobs} \sum_{i=1}^\numobs \wi f_j(\xij) \big| &
\leq \UniCon \; \big \{ \gamcrit^2 \; \|f_j\|_\Hil + \gamcrit \;
\|f_j\|_\numobs \big \} \qquad \mbox{for all $f_j \in \Hil$}
\end{align}
with probability greater than $1 - \plaincon_1 \exp(-\plaincon_2
\numobs \gamcrit^2)$.  By combining this result with our choice of
$\gamcrit$ and the union bound, the claimed bound then follows on
$\mprob[\TailEvent(\gamcrit)]$.

If $f_j = 0$, then the claim~\eqref{EqnNetherlandsRandom} is trivial.
Otherwise we renormalize $f_j$ by defining $g_j \defn
f_j/\|f_j\|_\Hil$, and we write
\begin{align*}
\frac{1}{\numobs} \sum_{i=1}^\numobs \wi f_j(\xij) & = \|f_j\|_\Hil \;
\frac{1}{\numobs} \sum_{i=1}^\numobs \wi g_j(\xij) \; \leq \;
\|f_j\|_\Hil \; \RvarHat \big(\w; \|g_j\|_\numobs, \Hil \big),
\end{align*}
where the final inequality uses the definition~\eqref{EqnDefnRvar},
and the fact that $\|g_j\|_\Hil = 1$.  We now split the analysis into
two cases: (1) $\|g_j\|_\numobs \leq \gamcrit$, and (2)
$\|g_j\|_\numobs > \gamcrit$.

\paragraph{Case 1: $\|g_j\|_\numobs \leq \gamcrit$.}  In this case,
it suffices to upper bound the quantity $\RvarHat (\w;
\gamma_\numobs, \Hil)$.  Note that $\|g_j\|_\Hil = 1$ and recall
definition~\eqref{EqnDefnRvar} of the random variable $\RvarHat$.
On one hand, since $\gamcrit \geq \delhatuni{j}$ by
Lemma~\ref{LemKerEquiv}, the definition of $\delhatuni{j}$ implies
that $\LocGaussEmp(\gamcrit, \Hil) \; \leq 4 \,\gamcrit^2$, and
hence
\begin{align*}
  \Exs[\RvarHat(w; \gamcrit; \Hil)] = \LocGaussEmp(\gamcrit; \Hil)
  & \leq 4 \gamcrit^2.
\end{align*}
Applying the bound~\eqref{EqnLedOne} from Lemma~\ref{LemGconc} with
$\delta = \gamcrit = t$, we conclude that $\RvarHat(w; \gamcrit;
\Hil) \leq \UniCon \; \gamcrit^2$ with probability at least $1 -
\plaincon_1 \exp \big \{-\plaincon_2 \numobs \gamcrit^2 \big \}$,
which completes the proof in the case where \mbox{$\|g\|_n \leq
  \gamcrit$.}

\paragraph{Case 2: $\|g_j\|_\numobs > \gamcrit$.}
In this case, we study the random variable $\RvarHat(w; r_j; \Hil)$
for some $r_j > \gamcrit$.  Our intermediate goal is to prove the
bound
\begin{align}
\label{EqnInterTailB}
\mprob \biggr[\RvarHat(w; r_j; \Hil) \geq \UniCon \, r_j \, \gamcrit
  \biggr] & \leq \plaincon_1 \exp \big \{-\plaincon_2 \numobs
\gamcrit^2 \big \}.
\end{align}
Applying the bound~\eqref{EqnLedOne} with $t = r_j$ and $\delta =
\gamcrit$, we are guaranteed an upper bound of the form
\mbox{$\RvarHat(w; r_j; \Hil) \leq \LocGaussEmp(r_j, \Hil) + r_j
  \, \gamcrit$} with probability at least $1-\plaincon_1 \exp \big
(-\plaincon_2 \numobs \gamcrit^2)$.  In order to complete the proof,
we need to show that $\LocGaussEmp(r_j, \Hil) \leq r_j \: \gamcrit$.
Since $r_j > \gamcrit > \delhatuni{j}$, we have
\begin{align*}
\LocGaussEmp(r_j, \Hil) & = \frac{r_j}{\delhatuni{j}} \Exs_{\gvar}
\big[ \sup_{\substack{\|g_j\|_\numobs \leq \delhatuni{j}
      \\ \|g_j\|_\Hil \leq \frac{\delhatuni{j}}{r_j}}} \big |
  \frac{1}{\numobs} \sum_{i=1}^\numobs \gvar_i g_j(\xij) \big| \big]
\; \leq \; \frac{r_j}{\delhatuni{j}} \LocGaussEmp(\delhatuni{j}, \Hil)
\; \leq 4 \, r_j \delhatuni{j},
\end{align*}
where the final inequality uses the fact that
$\LocGaussEmp(\delhatuni{j}, \Hil)\leq 4 \, \delunihat{j}^2$.  On the
event $\AnnoyEvent(\gamcrit)$ from Lemma~\ref{LemKerEquiv}, we have
$\delhatuni{j} \leq \gamcrit$, from which the
claim~\eqref{EqnInterTailB} follows.

We now use the bound~\eqref{EqnInterTailB} to prove the
bound~\eqref{EqnNetherlandsRandom}, in particular via a ``peeling''
operation over all choices of $r_j = \|f_j\|_\numobs/\|f_j\|_\Hil$.
(See van de Geer~\cite{vandeGeer} for more details on these peeling
arguments.)  We claim that it suffices to consider $r_j\leq 1$.  It is
equivalent to show that $\|g_j\|_\numobs \leq 1$ for any $g_j \in
\Ball_\Hil(1)$.  Since $\|g_j\|_\infty \leq \|g_j\|_{\Hil} \leq 1$, we
have $\|g_j\|_\numobs^2 = \frac{1}{\numobs} \sum_{i=1}^\numobs
g_j^2(x_{ij}) \; \leq 1$, as required.  Now define the event
\begin{align*}
\TailEvent_j(\gamcrit) & \defn \biggr \{ \exists f_j \in \Ball_\Hil(1)
\, \mid \, \big| \frac{1}{\numobs} \sum_{i=1}^\numobs \wi f_j(\xij)
\big| > 8 \; \|f_j\|_\Hil \; \gamcrit \;
\frac{\|f_j\|_\numobs}{\|f_j\|_\Hil}, \mbox{ and }
\frac{\|f_j\|_\numobs}{\|f_j\|_\Hil} \in (\gamcrit, 1] \biggr \}.
\end{align*}
and the sets $S_m \defn \big \{ 2^{m-1} \gamcrit \leq
\frac{\|f_j\|_\numobs}{\|f_j\|_\Hil} \leq 2^m \gamcrit \big \}$ for $m
= 1, 2, \ldots, M$.  By choosing $M = 2\log_2(1/\gamcrit)$, we ensure
that $2^M \gamcrit \geq 1$, and hence that if the event
$\TailEvent_j(\gamcrit)$ occurs, then it must occur for function $f_j$
belonging to some $S_m$, so that we have a function $f_j$ such that
$\frac{\|f_j\|_\numobs}{\|f_j\|_\Hil} \leq t_m \defn 2^{m} \gamcrit$,
and
\begin{align*}
\big| \frac{1}{\numobs} \sum_{i=1}^\numobs \wi f_j(\xij) \big| & > \,
8 \; \|f_j\|_\Hil \: \gamcrit \;
\frac{\|f_j\|_\numobs}{\|f_j\|_\Hil} \; \geq \UniCon \; \|f_j\|_\Hil
\: t_m,
\end{align*}
which implies that $\RvarHat(w; t_m, \Hil) \geq 4
t_m$. Consequently, by union bound and the tail
bound~\eqref{EqnInterTailB}, we have
\begin{align*}
\mprob[\TailEvent_j(\gamcrit)] & \leq M \; \plaincon_1 \exp \big
\{-\plaincon_2 \numobs \gamcrit^2 \big \} \; \leq \; \plaincon_1 \exp
\big \{ -\plaincon'_2 \numobs \gamcrit^2 \big \}
\end{align*}
by the condition $\numobs \gamcrit^2 = \Omega(\log(1/\gamcrit))$,
which completes the proof.


\section{Proof of Lemma~\ref{LemConeBound}}
\label{AppLemConeBound}

Define the function
\begin{align*}
\LossTil(\Delta) & \defn \frac{1}{2 \numobs} \sum_{i=1}^\numobs
\big(\wi + \mu +\ybar -\Delta(\exi) \big)^2 + \regpar \Lpeone{\fstar +
  \Delta} + \regpartwo \Hilone{\fstar + \Delta}
\end{align*}
and note that by definition of our $M$-estimator, the error function
$\DeltaHat \defn \fhat - \fstar$ minimizes $\LossTil$.  From the
inequality $\LossTil(\DeltaHat) \leq \LossTil(0)$, we obtain the upper
bound $\frac{1}{2}\|\DeltaHat\|_\numobs^2 \leq T_1 + T_2$, where
\begin{align*}
T_1 & \defn \big | \frac{1}{\numobs} \sum_{i=1}^\numobs
{\wi\DeltaHat(\exi)} \big | + |\ybar-\mu | \big | \frac{1}{\numobs}
\sum_{i=1}^\numobs {\DeltaHat(\exi)} \big |, \quad \mbox{and} \\
T_2 & \defn \regpar \sum_{j=1}^\pdim \big \{ \|\fstar_j \|_\numobs -
\|\fstar_j + \DeltaHat_j\|_\numobs \big\} + \regpartwo
\sum_{j=1}^\pdim \big \{ \|\fstar_j\|_\Hil - \|\fstar_j +
\DeltaHat_j\|_\Hil \big \}.
\end{align*}
Conditioned on the event $\mathcal{C}(\gamma_n)$, we have the bound
$|\ybar-\mu | \big | \frac{1}{\numobs} \sum_{i=1}^\numobs
{\DeltaHat(\exi)} \big | \leq \sqrt{\s} \gamcrit
\|\DeltaHat\|_\numobs$, and hence $\frac{1}{2}\|\DeltaHat\|_\numobs^2
\leq T_2 + \big | \frac{1}{\numobs} \sum_{i=1}^\numobs
     {\wi\DeltaHat(\exi)} \big | + \sqrt{\s} \gamcrit
     \|\DeltaHat\|_\numobs$, or equivalently
\begin{align}
\label{EqnParis}
0 \; \leq \; \frac{1}{2} \big( \|\DeltaHat\|_\numobs - \sqrt{\s}
\gamcrit \big)^2 & \leq T_2 + \big | \frac{1}{\numobs}
\sum_{i=1}^\numobs {\wi\DeltaHat(\exi)} \big | + \frac{1}{2} \s
\gamcrit^2.
\end{align}

It remains to control the term $T_2$.  On one hand, for any $j \in
\Sbar$, we have
\begin{equation*}
\|\fstar_j \|_\numobs - \|\fstar_j + \DeltaHat_j\|_\numobs \; = \;
-\|\DeltaHat_j\|_\numobs, \quad \mbox{and} \quad \|\fstar_j \|_\Hil -
\|\fstar_j + \DeltaHat_j\|_\Hil \; = \; -\|\DeltaHat_j\|_\Hil.
\end{equation*}
On the other hand, for any $j \in \Sset$, the triangle inequality
yields $\|\fstar_j \|_\numobs - \|\fstar_j + \DeltaHat_j\|_\numobs
\leq \|\DeltaHat_j\|_\numobs$, with a similar inequality for the terms
involving $\|\cdot\|_\Hil$.  Combined with the bound~\eqref{EqnParis},
we conclude that
\begin{align}
\label{EqnFriday}
0 & \leq \frac{1}{\numobs} \sum_{i=1}^\numobs \wi \DeltaHat(\exi)
+ \regpar \big \{ \Lpeone{\DeltaHat_\Sset} - \Lpeone{\DeltaHat_\Sbar}
\big\} + \regpartwo \big \{ \Hilone{\DeltaHat_\Sset} -
\Hilone{\DeltaHat_\Sbar} \big \} + \frac{1}{2}\s \gamcrit^2.
\end{align}
Recalling our conditioning on the event $\TailEvent(\gamcrit)$, by Lemma~\ref{LemNetherlands}, we have the upper bound
\begin{align*}
\big| \frac{1}{\numobs} \sum_{i=1}^\numobs \wi \DeltaHat(\exi)| &
\leq 8 \; \big \{ \gamcrit \Lpeone{\DeltaHat} +
\gamcrit^2 \Hilone{\DeltaHat} \big \}.
\end{align*}
Combining with the inequality~\eqref{EqnFriday} yields
\begin{align*}
0 & \leq 8 \; \big \{ \gamcrit \Lpeone{\DeltaHat} +
\gamcrit^2 \Hilone{\DeltaHat} \big \} + \regpar \big \{
\Lpeone{\DeltaHat_\Sset} - \Lpeone{\DeltaHat_\Sbar} \big \} +
\regpartwo \big \{ \Hilone{\DeltaHat_\Sset} - \Hilone{\DeltaHat_\Sbar}
\big \} + \frac{1}{2}\s \gamcrit^2\\
& \leq \frac{\regpar}{2} \Lpeone{\DeltaHat} + \frac{\regpartwo}{2}
\Hilone{\DeltaHat} + \regpar \big \{ \Lpeone{\DeltaHat_\Sset} -
\Lpeone{\DeltaHat_\Sbar} \big \} + \regpartwo \big \{
\Hilone{\DeltaHat_\Sset} - \Hilone{\DeltaHat_\Sbar} \big \} + \frac{1}{2}\s \gamcrit^2,
\end{align*}
where we have recalled our choices of $(\regpar, \regpartwo)$.
Finally, re-arranging terms yields the claim~\eqref{EqnConeBound}.


\section{Proof of Lemma~\ref{LemLowerBound}}
\label{AppLemLowerBound}

Recalling the definition~\eqref{EqnSuperGclass} of the function class
$\Gclass(\regpar, \regpartwo)$ and the critical radius $\tildelcrit$
from equation~\eqref{EqnDefnTildecrit}, we define the function class
\mbox{$\Gclass'(\regpar, \regpartwo, \tildelcrit) \defn \big \{ h \in
  \Gclass(\regpar, \regpartwo)\;|\; \|h\|_2 = \tildelcrit \big \}$,}
and the alternative event
\begin{align*}
\Bevent'(\regpar, \regpartwo) & \defn \big \{ \{ \|h\|_\numobs^2 \geq
\tildelcrit^2/2 \quad \mbox{for all $h \in \Gclass'(\regpar,
  \regpartwo, \tildelcrit)$} \big \}.
\end{align*}
We claim that it suffices to show that $\Bevent'(\regpar, \regpartwo)$
holds with probability at least \mbox{$1 - \plaincon_1 \exp(-
  \plaincon_2 \numobs \gamcrit^2)$.} Indeed, given an arbitrary
non-zero function $g \in \Gclass(\regpar, \regpartwo)$, consider the
rescaled function $h = \frac{\tildelcrit}{\|g\|_2} g$.  Since $g \in
\Gclass(\regpar, \regpartwo)$ and $\Gclass(\regpar, \regpartwo)$ is
star-shaped, we have $h \in \Gclass(\regpar, \regpartwo)$, and also
$\|h\|_2 = \tildelcrit$ by construction.  Consequently, when the event
$\Bevent'(\regpar, \regpartwo)$ holds, we have $\|h\|_\numobs^2 \geq
\tildelcrit^2/2$, or equivalently $\|g\|_\numobs^2 \geq \|g\|_2^2/2$,
showing that $\Bevent(\regpar, \regpartwo)$ holds. Accordingly, the
remainder of the proof is devoted to showing that $\Bevent'(\regpar,
\regpartwo)$ holds with probability greater than $1 - \plaincon_1
\exp(-\plaincon_2 \numobs \gamcrit^2)$.  Alternatively, if we define
the random variable $Z_\numobs(\Gclass') \defn \sup_{f \in
  \Gclass'}\big \{\tildelcrit^2 - \frac{1}{\numobs} \sum_{i=1}^\numobs
{f^2(x_i)} \big \}$, then it suffices to show that $Z_n(\Gclass') \leq
\tildelcrit^2/2$ with high probability.

Recall from Section~\ref{SecProofThmLower} the definition of a
covering set; here we use the notion of a proper covering, which
restricts the covering to use only members of the set $\GenSet$.
Letting $\Nprop(\epsilon; \GenSet, \rho)$ denote the propert covering
number, it can be shown that $\Nprop(\epsilon; \GenSet, \rho) \leq
N(\epsilon; \GenSet, \rho) \leq \Nprop(\epsilon/2; \GenSet, \rho)$.
Now let $g^1, \ldots, g^N$ be a minimal $\tildelcrit/8$-proper
covering of $\Gclass'$ in the $L^2(\mathbb{P}_n)$-norm, so that for
all $f \in \Gclass'$, there exists $g=g^k \in \Gclass'$ such that
$\|f-g\|_\numobs \leq \tildelcrit/8$.  We can then write
\begin{equation*}
\tildelcrit^2 - \frac{1}{n}\sum_{i=1}^{n}{f^2(x_i)} = \big \{
\tildelcrit^2 - \frac{1}{n}\sum_{i=1}^{n}{g^2(x_i)}\big \} + \big \{
\frac{1}{n}\sum_{i=1}^{n}{(g^2(x_i) - f^2(x_i))}\big \}.
\end{equation*}	
By the Cauchy-Schwartz inequality, we have
\begin{align*}
\frac{1}{\numobs} \sum_{i=1}^\numobs {(g^2(x_i) - f^2(x_i))} & =
\frac{1}{\numobs} \sum_{i=1}^\numobs {(g(x_i) - f(x_i))(g(x_i) +
  f(x_i))} \\
& \leq \sqrt{\frac{1}{\numobs} \sum_{i=1}^\numobs {(g(x_i) -
    f(x_i))^2}} \sqrt{\frac{1}{\numobs} \sum_{i=1}^\numobs (f(x_i) +
  g(x_i))^2} \\
& = \| g - f\|_\numobs \; \sqrt{\frac{1}{n}\sum_{i=1}^{n}{(f(x_i) +
    g(x_i))^2}}.
\end{align*}
By our choice of the covering, we have $\|g - f \|_\numobs \leq
\tildelcrit/8$.  On the other hand, we have
\begin{align*}
\sqrt{\frac{1}{\numobs} \sum_{i=1}^\numobs (f(x_i) + g(x_i))^2} & \leq
\sqrt{2 \|f\|_\numobs^2 + 2 \|g\|_\numobs^2} \; \leq \; \sqrt{4
  \tildelcrit^2} \; = \; 2 \tildelcrit,
\end{align*}
where the final inequality follows since $\|f\|_\numobs =
\|g\|_\numobs = \tildelcrit$.  Overall, we have established the upper
bound $\frac{1}{\numobs}\sum_{i=1}^\numobs {(g^2(x_i) - f^2(x_i))}
\leq \frac{\tildelcrit^2}{4}$, and hence shown that
\begin{align*}
Z_n(\Gclass') & \leq \max_{g^1, g^2, \ldots, g^N} \big \{
\tildelcrit^2 - \frac{1}{n}\sum_{i=1}^{n}{(g^k(x_i))}\big \} +
\frac{\tildelcrit^2}{4},
\end{align*}
where $N = \Nprop(\tildelcrit/8,\Gclass',\| \cdot\|_\numobs)$.  For
any $g$ in our covering set, since $g^2(x_i) \geq 0$, we may apply
Theorem 3.5 from Chung and Lu~\cite{Chung06} with $t =
\tildelcrit^2/4$ to obtain the one-sided tail bound
\begin{equation}
\label{EqnOneSide}
\mprob[\tildelcrit^2 - \frac{1}{\numobs} \sum_{i=1}^\numobs g^2(x_i)
  \geq \frac{\tildelcrit^2}{4}] \leq \exp \big(- \frac{\numobs
  \tildelcrit^4}{32 \Exs[g^4(x)]} \big ),
\end{equation}	
where we used the upper bound $\mbox{var}(g^2(x)) \leq
\mathbb{E}[g^4(x)]$. Next using the fact that the variables $\{
g_j(x_j)\}_{j=1}^\pdim$ are independent and zero-mean, we have
\begin{eqnarray*}
\Exs[g^4(x)] & = & \sum_{j=1}^\pdim \Exs[g_j^4(x_j)] + {4 \choose 2}
\sum_{j \neq k} \Exs[[g_j^2(x_j)] \Exs[g_k^2(x_k)] \\
& \leq & 4 \sum_{j=1}^\pdim {\mathbb{E}[g_j^2(x_j)]} + 6
\sum_{j=1}^\pdim {\Exs[g_j^2(x_j)]} \sum_{k=1}^\pdim
    {\Exs[g_k^2(x_k)]} \\
& \leq & 4 \tildelcrit^2 + 6 \tildelcrit^4\\
& \leq & 10 \tildelcrit^2,
\end{eqnarray*}	
where the second inequality follows since $\|g_j\|_{\infty} \leq
\|g_j\|_{\Hil} \leq 2$ for each $j$.  Combining this upper bound on
$\Exs[g^4(x)]$ with the earlier tail bound~\eqref{EqnOneSide} and
applying union bound yields
\begin{equation}
\label{EqnProb}
\mprob[ \max_{k=1,2,...,N} \big \{ \tildelcrit^2 - \frac{1}{\numobs}
  \sum_{i=1}^\numobs g^2(x_i) \big \} \geq \frac{\tildelcrit^2}{4}]
\leq \exp \big(\log \Nprop (\tildelcrit/8, \Gclass', \| \cdot
\|_\numobs) - \frac{\numobs \tildelcrit^2}{320} \big).
\end{equation}

It remains to bound the covering entropy $\log \Nprop(\tildelcrit/8,
\Gclass', \|\cdot\|_\numobs)$.  Since the proper covering entropy
$\log \Nprop(\tildelcrit/8, \Gclass', \|\cdot\|_\numobs)$ is at most
\mbox{$\log N(\tildelcrit/16, \Gclass', \| \cdot\|_\numobs)$,} it
suffices to upper bound the usual covering entropy.  Viewing the
samples $(x_1,x_2,...,x_n)$ as fixed, let us define the zero-mean
Gaussian process $\{W_g, g \in \Gclass'\}$ via $W_g \defn
\frac{1}{\sqrt{\numobs}} \sum_{i=1}^\numobs \newgvar_i g(x_i)$, where
the variables $\{\newgvar_i\}_{i=1}^\numobs$ are i.i.d. standard
Gaussian variates. By construction, we have $\mbox{var}[(W_g -W_f))] =
\|g-f\|_\numobs^2$. Consequently, by the Sudakov
minoration~\cite{Pisier89}, for all $\epsilon > 0$, we have $\epsilon
\sqrt{\log N(\epsilon; \Gclass', \|\cdot\|_\numobs)} \leq 4
\Exs_\newgvar[\sup_{g \in \Gclass'} W_g]$.  Setting $\epsilon =
\tildelcrit/16$ and performing some algebra, we obtain the upper bound
\begin{equation}
\label{EqnSudakov}
\frac{1}{\sqrt{\numobs}}\sqrt{\log N(\tildelcrit/16; \Gclass',
  \|\cdot\|_\numobs)} \leq \frac{64}{\tildelcrit} \mathbb{E}_\newgvar
     [ \sup_{g \in \Gclass'} \frac{1}{\numobs} \sum_{i=1}^\numobs
       \newgvar_i g(x_i)].
\end{equation}

The final step is to upper bound the Gaussian complexity $
\Exs_\newgvar[\sup \limits_{g \in \Gclass'}\frac{1}{\numobs}
  \sum_{i=1}^\numobs \newgvar_i g(x_i)]$.  In the proof of
Lemma~\ref{LemNetherlands}, we showed that for any co-ordinate $j \in
\{1, 2, \ldots, \pdim \}$, the univariate Gaussian complexity is upper
bounded as
\begin{align*}
\Exs \big[ \sup_{\substack{\|g_j\|_\numobs \leq r_j \\ \|g_j\|_\Hil
      \leq R_j}} \frac{1}{\numobs} \sum_{i=1}^\numobs \newgvar_i
  g_j(x_{ij}) \big] \; \leq \; C \, \big \{ \gamcrit \, r_j +
\gamcrit^2 R_j \big \}.
\end{align*}
Summing across co-ordinates and recalling the fact that the constant
$C$ may change from line to line, we obtain the upper bound
\begin{align*}
\Exs_\newgvar[\sup \limits_{g \in \Gclass'}\frac{1}{\numobs}
  \sum_{i=1}^\numobs \newgvar_i g(x_i)] & \leq
 C \sup_{g \in \Gclass'}
\big \{ \gamcrit \; \|g\|_{1, \numobs} + \gamcrit^2 \|g\|_{1, \Hil}
\big \} \\
& \stackrel{(a)}{\leq} C \sup_{g \in \Gclass'} \big \{ 4 \gamcrit \;
\|g_\Sset\|_{1, \numobs} + 4 \gamcrit^2 \|g_\Sset\|_{1, \Hil} +
\frac{1}{32} \s \regpartwo \big \} \\
& \stackrel{(b)}{\leq} C \sup_{g \in \Gclass'} \big \{ \gamcrit \;
\|g_\Sset\|_{1, \numobs} + \s \regpartwo \big \} \\
& \stackrel{(c)}{\leq} C \sup_{g \in \Gclass'} \biggr \{ \gamcrit \; [2
    \sqrt{\s} \|g\|_2 + \s \gamcrit ] + \s \regpartwo \biggr \},
\end{align*}
where step (a) uses inequality~\eqref{EqnGclassCone} in the definition
of $\Gclass'$; step (b) uses the inequality $\|g_j\|_\Hil \leq 2$ for
each co-ordinate and hence $\|g_\Sset\|_{1, \Hil} \leq 2 \s$, and our
choice of regularization parameter $\regpartwo \geq \gamcrit^2$; and
step (c) uses inequality~\eqref{EqnGclassSet} in the definition of
$\Gclass'$.  Since $\|g\|_2 = \tildelcrit$ for all $g \in \Gclass'$,
we have shown that
\begin{align}
\label{EqnBlueJay}
\Exs_\newgvar[\sup \limits_{g \in \Gclass'}\frac{1}{\numobs}
  \sum_{i=1}^\numobs \newgvar_i g(x_i)] & \leq C \big \{ \s \gamcrit^2
+ \sqrt{\s} \gamcrit \tildelcrit + \s \regpartwo \big\} \;
\stackrel{(d)}{\leq} \; C \big \{ \frac{\tildelcrit^2}{B^2} +
\frac{\tildelcrit^2}{B} \big \},
\end{align}
where inequality (d) follows from our choice~\eqref{EqnDefnTildecrit}
of $\tildelcrit$, and the constant $B$ can be chosen as large as we
please.  In particular, by choosing $B$ sufficiently large, and
combining the bound~\eqref{EqnBlueJay} with the Sudakov
bound~\eqref{EqnSudakov}, we can ensure that
\begin{align*}
\frac{1}{\numobs} \; \log N(\tildelcrit/16; \Gclass',
\|\cdot\|_\numobs) & \leq \frac{\tildelcrit^2}{640}.
\end{align*}
Combined with the earlier tail bound~\eqref{EqnProb}, we conclude
that
\begin{equation*}
\mprob [ \max_{k=1,2,...,N} \big \{ \tildelcrit^2 - \frac{1}{\numobs}
  \sum_{i=1}^\numobs g^2(x_i) \big \} \geq \frac{\tildelcrit^2}{4}]
\leq \exp \big(- \frac{ \numobs \tildelcrit^2}{640} \big),
\end{equation*}
which completes the proof of Lemma~\ref{LemLowerBound}.


\section{Proof of Lemma~\ref{LemMetEnt}}
\label{AppLemMetEnt}

\paragraph{Proof of part (a):}
Let $\SpecM = M(\frac{\delta}{\sqrt{\spindex}}; \Ball_\Hil(1),
\|\cdot\|_2) - 1$, and define $\IndSet = \{0,1, \ldots, \SpecM \}$.
Using \mbox{$\|u\|_0 = \sum_{j=1}^\pdim \Ind[u_j \neq 0]$} to denote
the number of non-zero components in a vector, consider the set
\begin{align}
\label{EqnDefnHset}
\Hset & \defn \big \{ u \in \IndSet^\pdim \, \mid \, \|u\|_0 =
\spindex \big \}.
\end{align}
Note that this set has cardinality $|\Hset| = {\pdim \choose \spindex}
\SpecM^\spindex$, since any element is defined by first choosing $\s$
co-ordinates are non-zero, and then for each co-ordinate, choosing
non-zero entry from a total of $\SpecM$ possible symbols.

For each $j = 1, \ldots, \pdim$, let $\{0, f_j^{1}, f_j^2, \ldots,
f_j^\SpecM \}$ be a $\delta/\sqrt{\spindex}$-packing of
$\Ball_\Hil(1)$.  Based on these packings of the univariate function
classes, we can use $\Hset$ to index a collection of functions
contained inside $\Fclass$.  In particular, any $u \in \Hset$ uniquely
defines a function \mbox{$g^u = \sum_{j=1}^d g_j^{u_j} \in \Fclass$,}
with elements
\begin{align}
g^{u_j}_j & = \begin{cases} f^{u_j}_j & \mbox{if $u_j \neq 0$} \\
                            0         & \mbox{otherwise.}
	      \end{cases}
\end{align}
Since $\|u\|_0 = \spindex$, we are guaranteed that at most $\spindex$
co-ordinates of $g$ are non-zero, so that $g \in \Fclass$.

  Now consider two functions $g^u$ and $h^v$ contained within the
class $\{g^u, u \in \Hset \}$.  By definition, we have
\begin{align}
\label{EqnPackLower}
\|g^u - h^v\|_2^2 & = \sum_{j=1}^\pdim \|f_j^{u_j} - f_j^{v_j}\|_2^2
\; \geq \; \frac{\delta^2}{\spindex} \sum_{j=1}^\pdim \Ind[u_j \neq
v_j],
\end{align}

Consequently, it suffices to establish the existence of a ``large''
subset $\Aset \subset \Hset$ such that the Hamming metric $\rho_H(u,
v) \defn \sum_{j=1}^\pdim \Ind[u_j \neq v_j]$ is at least $\spindex/2$
for all pairs $u, v \in \Aset$, in which case we are guaranteed that
$\|g - h\|^2_2 \geq \delta^2$.  For any $u \in \Hset$, we observe that
\begin{align*}
\biggr| \big \{ v \in \Hset \, \mid \, \rho_H(u,v) \leq \frac{s}{2}
\big \} \biggr| & \leq {\pdim \choose \frac{\spindex}{2}} \; (\SpecM +
1)^{\frac{\spindex}{2}}.
\end{align*}
This bound follows because we simply need to choose a subset of size
$\s/2$ where $u$ and $v$ agree, and the remaining $\s/2$ co-ordinates
can be chosen arbitrarily in $(\SpecM + 1)^{\frac{\spindex}{2}}$ ways.
For a given set $\Aset$, we write $\rho_H(u, \Aset) \leq
\frac{\spindex}{2}$ if there exists some $v \in \Aset$ such that
$\rho_H(u,v) \leq \frac{\spindex}{2}$.  Using this notation, we have
\begin{align*}
\biggr| \big \{ u \in \Hset \, \mid \, \rho_H(u, \Aset) \leq
\frac{\spindex}{2} \big \} \biggr| & \leq |\Aset| \; {\pdim \choose
\frac{\spindex}{2}} \; (\SpecM + 1)^{\frac{\spindex}{2}} \;
\stackrel{(a)}{<} |\Hset|,
\end{align*}
where inequality (a) follows as long as
\begin{align*}
|\Aset| & \leq N^* \; \defn \frac{1}{2} \frac{ {\pdim \choose
 \spindex}}{{\pdim \choose \frac{\spindex}{2}}} \; \frac{
 \SpecM^\spindex}{(\SpecM + 1)^{\spindex/2}}.
\end{align*}
Thus, as long as $|\Aset| \leq N^*$, there must exist some element $u
\in \Hset$ such that $\rho_H(u, \Aset) > \frac{\spindex}{2}$, in which
case we can form the augmented set $\Aset \cup \{u\}$.  Iterating this
procedure, we can form a set with $N^*$ elements such that
$\rho_H(u,v) \geq \frac{\spindex}{2}$ for all $u,v \in \Aset$.

Finally, we lower bound $N^*$.  We have
\begin{align*}
N^* & \stackrel{(i)}{\geq} \frac{1}{2} \, \big(\frac{\pdim -
\spindex}{\spindex/2} \big)^{\frac{\spindex}{2}} \;
\frac{(\SpecM)^\spindex}{(\SpecM +1)^{\spindex/2}} \\
& = \frac{1}{2} \, \big(\frac{\pdim - \spindex}{\spindex/2}
\big)^{\frac{\spindex}{2}} \; \SpecM^{\spindex/2} \big(\frac{
\SpecM}{\SpecM +1} \big)^{\spindex/2} \\
& \geq \frac{1}{2} \, \big(\frac{\pdim - \spindex}{\spindex/2}
\big)^{\frac{\spindex}{2}} \; \SpecM^{\spindex/2},
\end{align*}
where inequality (i) follows by elementary combinatorics (see Lemma 5
in the paper~\cite{RasWaiYu09} for details).  We conclude that for
$\spindex \leq \pdim/4$, we have
\begin{align*}
\log N^* & = \Omega \big( \spindex \log \frac{\pdim}{\spindex} +
\spindex \log M(\frac{\delta}{\sqrt{\spindex}}; \Ball_\Hil(1),
\|\cdot\|_2) \big),
\end{align*}
thereby completing the proof of Lemma~\ref{LemMetEnt}(a).

\paragraph{Proof of part (b):}
In order to prove part (b), we instead let $\SpecM = M(\frac{1}{2};
\Ball_\Hil(1), \|\cdot\|_2) - 1$, and then follow the same steps.
Since $\log \SpecM = \Omega(\mdim)$, we have the modified lower bound
\begin{align*}
\log N^* & = \Omega \big( \spindex \log \frac{\pdim}{\spindex} +
\spindex \mdim \big),
\end{align*}
Moreover, instead of the lower bound~\eqref{EqnPackLower}, we have
\begin{align}
\label{EqnPackLowerAlt}
\|g^u - h^v\|_2^2 & = \sum_{j=1}^\pdim \|f_j^{u_j} - f_j^{v_j}\|_2^2
\; \geq \; \frac{1}{4} \sum_{j=1}^\pdim \Ind[u_j \neq v_j] \; \geq \;
\frac{s}{8},
\end{align}
using our previous result on the Hamming separation.  Furthermore,
since $\|f_j\|_2 \leq \|f_j\|_\Hil$ for any univariate function, we
have the upper bound
\begin{align*}
\|g^u - h^v\|_2^2 & = \sum_{j=1}^\pdim \|f_j^{u_j} - f_j^{v_j}\|_2^2
\; \leq \; \sum_{j=1}^\pdim \|f_j^{u_j} - f_j^{v_j}\|_\Hil^2.
\end{align*}
By the definition~\eqref{EqnDefnHset} of $\Hset$, at most $2 \spindex$
of the terms $f_j^{u_j} - f_j^{v_j}$ can be non-zero.  Moreover, by
construction we have $\|f_j^{u_j} - f_j^{v_j}\|_\Hil \leq 2$, and hence
\begin{align*}
\|g^u - h^v\|_2^2 & \leq 8 \spindex.
\end{align*}
Finally, by rescaling the functions by $\sqrt{8} \, \delta/\sqrt{\s}$,
we obtain a class of $N^*$ rescaled functions $\{\gtil^u, u \in
\IndSet \}$ such that
\begin{equation*}
\|\gtil^u - \htil^v\|^2_2 \geq \delta^2, \quad \mbox{and} \quad
\|\gtil^u - \htil^v\|^2_2 \leq 64 \delta^2,
\end{equation*}
as claimed.

\section{Results for proof of Theorem~\ref{ThmFast}}
\label{AppThmFast}

The reader should recall from Section~\ref{SecKYCompare} the
definitions of the function classes $\MyNewClass$ and
$\MyNewSmallClass$. The function class $\MyNewSmallClass$ can be
parameterized by the two-dimensional sequence $(a_{j,k})_{j \in S\;, k
  \in \mathbb{N}}$ of co-efficients, and expressed in terms of
two-dimesnional sequence of basis functions $(\phi_{j,k})_{j \in S\;,
  k \in \mathbb{N}}$ and the sequence of eigenvalues $(\mu_k)_{k \in
  \mathbb{N}}$ for the univariate RKHS $\Hil$ as follows:
\begin{align*}
\MyNewSmallClass & \defn \big \{ f = \sum_{j \in \Sset}
\sum_{k=1}^\infty {a_{j,k}\phi_{j,k} } \, \mid \, \sum_{k=1}^\infty
    {\frac{a^2_{j,k}}{\mu_k}} \leq 1 \; \, \forall \; j \in \Sset\;
    \mbox{and}\;\|f\|_{\infty} \leq \GLOB \big \}.
\end{align*}
For any integer $M \geq 1$,  we also consider the truncated function
class
\begin{align*}
\TRUNCLASS & \defn \big \{ f = \sum_{j \in \Sset} \sum_{k=1}^M
           {a_{j,k}\phi_{j,k} } \, \mid \, \sum_{k=1}^\infty
           {\frac{a^2_{j,k}}{\mu_k}} \leq 1 \; \, \forall \; j \in
           \Sset\; \mbox{and}\;\|f\|_{\infty} \leq \GLOB \big \}.
\end{align*}
\blems
\label{LemSimple}
We have the inclusion $\TRUNCLASS \subseteq \big \{ f \in \MyHclass
\mid \, \sum_{j \in S}\sum_{k=1}^{M}|a_{j,k}| \leq B \, \sqrt{M}\big
\}$.  \elems
\spro 
Without loss of generality, let us assume that $S = \{1,2,...,\s \}$,
and consider a function $f = \sum_{j=1}^\s f_j \in \TRUNCLASS$.  Since
each $f_j$ acts on a different co-ordinate, we are guaranteed that
$\|f\|_\infty = \sum_{j = 1}^\kdim \|f_j\|_\infty$.  Consider any
univariate function $f_j = \sum_{k=1}^{M} a_{j,k} \phi_{j,k}$. We have
\begin{align*}
\sum_{k=1}^{M} |a_{j,k}| & \leq \sqrt{M} \, \biggr( \sum_{k=1}^M
a^2_{j,k} \biggr)^{1/2} \; \stackrel{(a)}{\leq} \; \sqrt{M} \; \big[
  \Exs[f_j^2(X_j)] \big]^{1/2} \; \leq \sqrt{M} \|f_j\|_\infty,
\end{align*}
where step (a) uses the fact that $\Exs[f_j^2(X_j)] =
\sum_{k=1}^\infty a^2_{j,k} \geq \sum_{k=1}^M a^2_{j,k}$ for any $M
\geq 1$.  Adding up the bounds over all co-ordinates, we obtain
\begin{align*}
\|a\|_1 \; = \; \sum_{j=1}^\kdim \sum_{k=1}^M |a_{j,k}| & \leq
\sqrt{M} \sum_{j=1}^\kdim \|f_j\|_\infty \; = \; \sqrt{M} \|f\|_\infty
\; \leq \; \sqrt{M} \GLOB,
\end{align*}
where the final step uses the uniform boundedness condition.  
\fpro

\subsection{Proof of Lemma~\ref{LemSquash}}
\label{AppLemRadGauss}

Recalling the definition of $\RvarHat(w; t, \HACK)$ stated
from~\eqref{GaussCompBound}, let us view it as a function of the
standard Gaussian random vector $(w_1, \ldots, w_\numobs)$.  It is
straightforward to verify that this variable is Lipschitz (with
respect to the Euclidean norm) with parameter at most
$t/\sqrt{\numobs}$.  Consequently, by concentration for Lipschitz
functions~\cite{Ledoux01}, we have
\begin{align*}
\mprob \big[\RvarHat(w; t, \HACK) \geq \Exs[\RvarHat(w; t,\HACK)] + 3
  t \UnifBoundRate \big ] & \leq \exp \big( - \frac{9 \numobs
  \UnifBoundRate^2}{2} \big).
\end{align*}
Next we prove an upper bound on the expectations
\begin{subequations}
\begin{align}
\LocGaussEmp(t; \HACK) & \defn \; \Exs_{w} \big[ \sup_{\substack{ g
      \in \HACK \\ \|g\|_\numobs \leq t}} \frac{1}{\numobs}
  \sum_{i=1}^\numobs w_i g(\exi) \big], \quad \mbox{and} \\
\LocGauss(t; \HACK) & \defn \; \Exs_{x,w} \big[ \sup_{\substack{ g
      \in \HACK \\ \|g\|_2 \leq t}} \frac{1}{\numobs}
  \sum_{i=1}^\numobs w_i g(\exi) \big].
\end{align}
\end{subequations}

\blems 
\label{LemSmallGauss}
Under the conditions of Theorem~\ref{ThmFast}, we have
\begin{eqnarray*}
\max \big \{ \LocGaussEmp(t; \HACK), \; \LocGauss(t; \HACK) \big \} &
\leq & 8 \GLOB \, C \,\sqrt{\frac{\s^{1/\alpha}\log \s}{\numobs}}.
\end{eqnarray*}
\elems

\spro
By definition, any function $g \in \HACK$ has support at most $2 \s$,
and without loss of generality (re-indexing as necessary), we assume
that $S = \{ 1,2,...,2 \s\}$.  We can thus view functions in $\HACK$
as having domain $\real^{2 \s}$, and we can an operator $\Phi$ that
maps from $\real^{2s}$ to $[\ell^2(\Nat)]^{2 s}$, via
\begin{equation*}
x \mapsto \Phi_{j,k}(x) = \phi_{j,k}(x_j), \qquad \mbox{for $j = 1,
  \ldots, 2 \s$, and $k \in \Nat$.}
\end{equation*}	
Any function in $g \in \HACK$ can be expressed in terms of
two-dimensional sequence $(a_{j,k})$ and the functions $(\Phi_{j,k})$
as $g(x) = g(x_1,x_2, \ldots ,x_{2 \s}) = \sum_{j=1}^{2 \s}
\sum_{k=1}^{\infty} \Phi_{j,k}(x) a_{j,k} \; = \;
\tracer{\Phi(x)}{a}$, where $\tracer{\cdot}{\cdot}$ is a convenient
shorthand for the inner product between the two arrays. \\

For any function $g \in \HACK$, triangle inequality yields the upper
bound
\begin{equation}
\label{EqnTri}
\sup_{g \in 2 \HACK} \frac{1}{n}|\sum_{i=1}^{n}{w_i
  \tracer{\Phi(x_i)}{a}}| \leq \underbrace{\sup_{g \in 2
    \HACK} \frac{1}{n}|\sum_{i=1}^{n}w_i
  \tracer{\Phi_{\cdot,1:M}(x_i)}{a_{\cdot, 1:M}}}_{\Aone} + \Atwo
\end{equation}
where $\Atwo \defn \sup_{g \in 2 \HACK}
\frac{1}{n}|\sum_{i=1}^{n} w_i \tracer{\Phi_{\cdot,M+1:\infty}(x_i)}{
  a_{\cdot,M+1:\infty}}|$.

\paragraph{Bounding term $\Exs_{x,w}[\Aone]$ and $\Exs_{w}[\Aone]$:} By H\"{o}lder's 
inequality and Lemma~\ref{LemSimple}, we have
\begin{align*}
\Aone & \leq \frac{1}{\sqrt{\numobs}} \sup_{g \in 2 \HACK}
\|a_{\cdot,1:M}\|_{1,1} \max_{j,k}
|\sum_{i=1}^{n}{\frac{\wi}{\sqrt{n}} \Phi_{j,k}(x_i)}| \; \leq \;
\frac{2 \;B \; \sqrt{M}}{\sqrt{\numobs}} \max_{j, k}
|\sum_{i=1}^{n}{\frac{\wi}{\sqrt{n}}\Phi_{j,k}(x_i)}|.
\end{align*}
By assumption, we have $|\Phi_{j,k}(x_i)| \leq C$ for all indices
$(i,j,k)$, implying that $\sum_{i=1}^{n}{\frac{\wi}{\sqrt{n}}
  \Phi_{j,k}(x_i)}$ is zero-mean with sub-Gaussian parameter bounded
by $C$ and we are taking the maximum of $2 \s \times M$ such
terms. Consequently, we conclude that
\begin{align}
\label{EqnAoneBound}
\Exs_{w}[\Aone] & \leq 8 \GLOB C \sqrt{\frac{M \log (M \s)}{\numobs}}.
\end{align}
The same bound holds for $\Exs_{x,w}[\Aone]$.

\paragraph{Bounding term $\Exs_{x,w}[\Atwo]$ and $\Exs_{w}[\Atwo]$:}
In order to control this term, we simply recognize that it corresponds
to the usual Gaussian complexity of the sum of $2 \s$ univariate
Hilbert spaces, each of which is an RKHS truncated to the
eigenfunctions $\{\mu_k\}_{k \geq M+1}$.
\begin{eqnarray*}  
\frac{1}{n}|\sum_{i=1}^{n} w_i \tracer{\Phi_{\cdot,M+1:\infty}(x_i)}{
  a_{\cdot,M+1:\infty}}| & \leq & \frac{1}{\sqrt{n}}\sum_{j=1}^{2
  \s}{|\sum_{k \geq M+1}{a_{j,k} \Phi_{j,k}(x)
    \sum_{i=1}^{n}{\frac{w_i}{\sqrt{n}}}} |}\\ & \leq &
\frac{C}{\sqrt{n}}\sum_{j=1}^{2 \s}{|\sum_{k \geq
    M+1}{\frac{a_{j,k}}{\sqrt{\mu_k}} \sqrt{\mu_k}
    \sum_{i=1}^{n}{\frac{w_i}{\sqrt{n}}}} |}\\ & \leq &
\frac{C}{\sqrt{n}}\sum_{j=1}^{2 \s}{\sqrt{\sum_{k \geq
      M+1}{\frac{a_{j,k}^2}{\mu_k}}} {\sqrt{\sum_{k \geq M+1}{\mu_k
        (\sum_{i=1}^{n}{\frac{w_i}{\sqrt{n}}})^2}}}},
\end{eqnarray*}
where the final inequality follows from Cauchy-Schwartz. Exploiting
the fact that $\sum_{k \geq M+1}{\frac{a_{j,k}^2}{\mu_k}} \leq 1$ for
all $j$, we have the bound
\begin{align}
\label{EqnAtwo}
\Exs_{w}[\Atwo] & \leq 4 C \s \; \frac{\sqrt{\sum_{k \geq M+1}{\mu_k}}}{\sqrt{n}}.
\end{align}
One again a similar bound holds for $\Exs_{x,w}[\Atwo]$.

Substituting the bound~\eqref{EqnAoneBound} and~\eqref{EqnAtwo} into
the upper bound~\eqref{EqnTri}, we conclude that
\begin{align*}
\LocGauss(2 \HACK) & \leq 4 B C \sqrt{\frac{M \log(M
    \s)}{n}} + 4 C \s \sqrt{\frac{\sum_{k \geq M+1}{\mu_k}}{n}} \\
& \leq 4 B C \sqrt{\frac{M \log(M \s)}{n}} + 4 C \s
\sqrt{\frac{M^{1-2 \alpha}}{n}},
\end{align*}
where the second inequality follows from the relation $\mu_k \simeq
k^{-2 \alpha}$.  Finally, setting $M = \s^{\frac{1}{\alpha}}$ yields
the claim.

Note that the same argument works for the
Rademacher complexity, since we only exploited the sub-Gaussianity of
the variables $w_i$.
\fpro

Returning to the proof of Lemma~\ref{LemSquash}, combining
Lemma~\ref{LemSmallGauss} with the bound ~\eqref{EqnLedOne} in
Lemma~\ref{LemGconc}:
\begin{align*}
\mprob \big[\RvarHat(w; t, \HACK) \geq 8 \GLOB \, C
  \,\sqrt{\frac{\s^{1/\alpha}\log \s}{n}} + 3 t \UnifBoundRate \big ]
& \leq \exp \big( - \frac{9 \numobs \UnifBoundRate^2}{2} \big).
\end{align*}
Since $\|g\|_\numobs \leq 2 \GLOB$ for any function $g \in \HACK$, the
proof Lemma~\ref{LemSquash} is completed using a peeling argument over
the radius, analogous to the proof of Lemma~\ref{LemNetherlands} (see
Appendix~\ref{AppLemNetherlands}).




\bibliographystyle{plain}
\bibliography{garvesh_dec11}

\end{document}